\documentclass[12pt]{amsart}
\usepackage{amscd,amsmath,amssymb,amsthm}
\newlength{\figboxwidth}
\newtheorem{thm}{Theorem}[section]
\newtheorem{cor}[thm]{Corollary}
\newtheorem{lem}[thm]{Lemma}
\newtheorem{prop}[thm]{Proposition}
\theoremstyle{definition}
\newtheorem{defn}[thm]{Definition}
\theoremstyle{remark}

\numberwithin{equation}{section}

\newcommand{\R}{\mathbb R}\newcommand{\C}{\mathbb C}\newcommand{\Z}{\mathbb Z}

\newcommand{\To}{\rightarrow}

\newcommand{\mc}[1]{{\sc #1}}
\newcommand{\G}{\Gamma}\newcommand{\g}{\gamma}
\newcommand{\fh}{\mathcal{H}}\newcommand{\cu}{\mathcal{U}}
\newcounter{bean}
\begin{document}

\title[Strengthening $(T)$]{Strengthening Kazhdan's Property $(T)$ by Bochner Methods}%
\author{David Fisher and Theron Hitchman}%
\email{fisherdm@indiana.edu, thitchma@williams.edu}%

\thanks{\noindent First author partially supported by NSF Grants DMS-0226121 and DMS-0541917.
Second author partially supported by NSF VIGRE grant DMS-0240058.}%

\date{\today}%
\begin{abstract}
In this paper, we propose a property which is a natural
generalization of Kazhdan's property $(T)$ and prove that many,
but not all, groups with property $(T)$ also have this property.

Let $\G$ be a finitely generated group. One definition of $\G$
having property $(T)$ is that $H^1(\G,\pi,\fh)=0$ where the
coefficient module $\fh$ is a Hilbert space and $\pi$ is a unitary
representation of $\G$ on $\fh$. Here we allow more general
coefficients and say that $\G$ has property $F \otimes H$ if
$H^1(\G,\pi_1{\otimes}\pi_2,F{\otimes}\fh)=0$ if $(F,\pi_1)$ is
any representation with $\dim(F)<\infty$ and $(\fh,\pi_2)$ is a
unitary representation.

The main result of this paper is that a uniform lattice in a
semisimple Lie group has property $F \otimes H$ if and only if it
has property $(T)$. The proof hinges on an extension of a
Bochner-type formula due to Matsushima-Murakami and Raghunathan.
We give a new and more transparent derivation of this formula as
the difference of two classical Weitzenb\"{o}ck formula's for two
different structures on the same bundle.  Our Bochner-type formula
is also used in our work on harmonic maps into continuum products
\cite{Fisher-Hitchman2,Fisher-Hitchman1}. Some further
applications of property $F \otimes H$ in the context of group
actions will be given in \cite{Fisher-Hitchman3}.
\end{abstract}
\maketitle
\section{Introduction and Statements of Results}
\label{section:intro}

Property $(T)$, introduced by Kazhdan in 1966 in \cite{Kazhdan-T},
plays a fundamental role in the study of discrete subgroups of Lie
groups and more general finitely generated groups.  In this paper
we propose a stronger property, which we call {\em property $F
\otimes H$}, which is a direct strengthening of one equivalent
definition of property $(T)$.  It follows from work of Delorme and
Guichardet that property $(T)$ is equivalent to the statement that
$H^1(\G, \pi, \fh)=0$ whenever $\fh$ is a Hilbert space and $\pi$
is a continuous unitary representation. This is easily seen to be
equivalent to the statement that any continuous affine isometric
$\G$ action on a Hilbert space has a fixed point, commonly called
property $FH$.

\begin{defn}
\label{definition:ft} Let $D$ be a topological group,  we say $D$
has property $F \otimes H$ if for every
\begin{enumerate}
\item $F$ a finite dimensional vector space over $\R$ or $\C$,
\item continuous homomorphism $\pi:D{\rightarrow}GL(F)$, \item
Hilbert space $\fh$, and \item continuous unitary representation
$\pi_2:D{\rightarrow}\cu(\fh)$,
\end{enumerate} the first continuous cohomology group
\noindent $H^1(\G,\pi_1{\otimes}\pi_2,F{\otimes}\fh)=0$.
\end{defn}

\noindent We refer to representations $\pi_1{\otimes}\pi_2$ as in
the definition above as {\em finite dimensional by unitary}
representations. 
Throughout this paper, cohomology for topological groups will be
continuous cohomology. It is easy to show that property $F \otimes
H$ is equivalent to having any continuous affine action of $D$ with
finite dimensional by unitary linear part has a fixed point. The
notation for property $F \otimes H$ evokes both the fixed point
property $FH$ generalized here and the tensor product in the finite
dimensional by unitary representations that appear in the
definition.

Throughout this paper a semisimple Lie group will be a connected
semisimple Lie group with finite center. We will not consider the
case of semisimple Lie groups which are not algebraic, see the end
of \S\ref{section:otherexamples} for discussion of related examples.
We will also always assume that $G$ is simply connected as an
algebraic group, i.e. that every cover of $G$ which is algebraic is
isomorphic to $G$. With this assumption, $G$ is the direct product
of it's simple factors \cite[Proposition I.1.4.10]{Margulis-Book}.
As it is obvious that property $F \otimes H$ passes to quotients,
there is no loss of generality in this assumption. The main result
of this paper is the following:

\begin{thm}
\label{theorem:ft} Let $D$ be a semisimple Lie group or a
cocompact lattice in such.  Then $D$ has property $F \otimes H$ if
and only if it has property $(T)$.
\end{thm}

\noindent{\bf On Theorem \ref{theorem:ft}:}
\begin{enumerate}

\item For $D$ a connected semisimple Lie group and $\pi_2$
irreducible, or more generally admissible, the desired vanishing
follows from known results on relative Lie algebra cohomology. It
does not appear possible to deduce the general case from this,
mainly because those results do not seem to yield estimates on the
size of solutions to cohomological equations, see \S
\ref{subsection:gkcohomologydirect}.

 \item For a general
locally compact group $D$ and a cocompact lattice $\Lambda<D$, it
follows from work of Blanc that if $D$ has $F \otimes H$ then
$\Lambda$ has $F \otimes H$ \cite{Blanc}. The converse (shown to
hold for property $(T)$ by Kazhdan) is not immediate.

\item We prove Theorem \ref{theorem:ft} for cocompact lattices and
prove it for the ambient group from that.  The proof does not
proceed by deducing $F \otimes H$ from $(T)$, but rather uses the
fact that $G$ and $\G$ have $(T)$ as long as $G$ has no factors
locally isomorphic to $SO(1,n)$ or $SU(1,n)$.

\item It is not true in general that property $F \otimes H$ is equivalent
to property $(T)$.  The easiest example is $SL(n,\Z){\ltimes}\Z^n$
for $n>2$, see \S\ref{section:otherexamples} below.

\item  $F \otimes H$ for cocompact lattices was announced in
\cite[Theorem 1.8]{Fisher-Hitchman1}.

\item The issue of whether non-uniform lattices have $F \otimes H$
seems quite subtle and will be addressed elsewhere.
\end{enumerate}

As mentioned above, the main step is proving Theorem
\ref{theorem:ft} for cocompact lattices.  We deduce $F \otimes H$
for $G$ from $F \otimes H$ for irreducible cocompact lattices in
$G{\times}G$. Our proof also depends on a use of the
Margulis-Corlette-Gromov-Schoen superrigidity theorems, which
classify all finite dimensional representations of the lattices we
consider \cite{Corlette-Annals,
Gromov-Schoen,Margulis-ICM,Margulis-Book}.

For cocompact lattices we translate the question into one
concerning the de Rham cohomology of the associated locally
symmetric space with coefficients in a flat $F{\otimes}\fh$
bundle. We then use a Bochner-Matsushima-type formula to prove an
estimate on the first eigenvalue of the Laplacian on one forms
that implies vanishing of first cohomology. In order to define a
Laplacian one requires a choice of metric on $F{\otimes}\fh$, and
our choice here is similar to the one made by Matsushima and
Murakami in \cite{Matsushima-Murakami}. In fact, the work here is
very close to the work in that paper, and we eventually reduce to
an estimate on eigenvalues of the \emph{same} finite dimensional
matrix as they do. These estimates are obtained by Raghunathan in
\cite{Raghunathan-1-vanishing}.  It is possible to derive the
formula we use by closely following the derivation in
\cite{Matsushima-Murakami}; however, we choose a different method.
From our point of view, the estimate follows by subtracting two
standard Weitzenb\"{o}ck formulas for different bundle structures
on the same vector bundle. The use of two bundle structures is
also a key technique in \cite{Matsushima-Murakami}, but the
relation of the computations there to standard differential
geometric computations is not clear. One benefit of our point of
view is that it makes it immediately clear why the negative Ricci
curvature of $M$ does not spoil the computations: it appears in
both formulas, and cancels when they are subtracted.  It is
worthwhile to compare our differential geometric interpretation of
the formula from \cite{Matsushima-Murakami} to a well known
differential geometric interpretation of an earlier Bochner-type
formula of Matsushima for trivial bundles
\cite{Matsushima-trivialrep}. This well-known differential
geometric interpretation is due to Calabi, was first described in
print in detail by Dodziuk \cite{Dodziuk}, and played a key role
in the development of geometric superrigidity
\cite{Jost-Yau,Mok-Siu-Yeung}.  As in those works, our Bochner
formula will be useful in non-linear settings, see
\cite{Fisher-Hitchman2,Fisher-Hitchman1}.

In \S\ref{subsection:higherdegrees}, we deduce some vanishing
results for $H^i(\G,\pi_1{\otimes}\pi_2 , F{\otimes}\fh)$ for more
restricted choices of $F$.  These results are probably not
optimal. A reader familiar with the literature on cohomology of
Lie groups and lattices might wonder why we do not pursue an
analysis based on relative Lie algebra cohomology and methods
closer to those of \cite{Borel-Wallach, Schmid-Advances,
Zuckerman-Annals}.  Approaches of this kind, whether directly to
show vanishing for $G$ cohomology, or indirectly to show vanishing
of $\Gamma$ cohomology via Matsushima-type formulas, encounter
analytic difficulties. We discuss the issues that arise in \S
\ref{section:gkcohomology}.  We remark that our work gives a new
interpretation of the Matsushima-Murakami formula for the
Laplacian on vector bundles of the form
$(K{\backslash}G{\times}F)\G$ which is used in the proof of
Matsushima's formula. Our more geometric interpretation of this
formula may be of interest even to those primarily interested in
traditional applications of Matsushima's formula.

In addition, the estimate on the Laplacian proved here using the
Bochner method is also used in our work on harmonic maps into
continuum products \cite{Fisher-Hitchman2, Fisher-Hitchman1}.

In \cite{Fisher-Hitchman3}, using results in \cite{Fisher-H1} we
prove a local rigidity theorem for many geometric actions of
groups with property $F \otimes H$ on compact manifolds. This can
be seen as a generalization of the theorem of the first author and
Margulis on local rigidity of isometric actions of groups with
property $(T)$ \cite{Fisher-Margulis1}.

Parts of this work were supported by visits to the Graduate Center
of the City University of New York, Indiana University and Rice
University. We thank these institutions for their hospitality and
support.  We would also like to thank Nicolas Monod for suggesting
the name property $F{\otimes}H$.

\section{Other groups \mc{and} property $F \otimes H$}
\label{section:otherexamples}

Other than the groups that arise in Theorem \ref{theorem:ft}, other
groups known to have property $(T)$ fall into a few distinct
classes. In this section we discuss which of these have or do not
have property $F \otimes H$.

Let $H$ be a Lie group with Levi decomposition $H=L{\ltimes}\R^n$
and assume that $L$ is semisimple with no compact factors and
property $(T)$ and that the $L$ representation on $\R^n$ has no
invariant vectors.  Then $H$ is known to have property $(T)$, see
e.g. \cite{delaHarpe-Valette-Book,Wang,Zimmer-Book}.  In many
cases, it is obvious that $H$ contains lattices, e.g.
$SL(n,\Z){\ltimes}\Z^n<SL(n,\R){\ltimes}\R^n$.  For a more
detailed discussion of $H$'s containing lattices, see
\cite{Valette-Pairpaper}. In any case, we have:

\begin{prop}
If $H$ is as in the preceding paragraph, then both $H$ and any
lattice in $H$ fail to have property $F \otimes H$.
\end{prop}

\begin{proof}
It is clear that $H$ admits an action on $\R^n$ with no fixed
points, by letting $(l,v)\cdot w= lw+v$ where $l\in L$ acts on $w$
by the representation defining the semidirect product. It is also
clear that this action has no fixed points when restricted to a
lattice in $H$. This can be viewed as providing non-trivial
cohomology classes in $H^1(H,\R^n)$ or $H^1(\G,\R^n)$.
\end{proof}

This observation extends to more general semidirect products
$L{\ltimes}U$ with $U$ unipotent, we leave the construction of more
examples to the reader. For all other groups we know of with
property $(T)$, the key fact seems to be the following simple lemma.

\begin{lem}
\label{lemma:stupid} Let $D$ be a group with property $(T)$ such
that every continuous finite dimensional representation of $D$ over
$\R$ or $\C$ is unitary. Then $D$ has property $F \otimes H$.
\end{lem}

The lemma allows us to deduce that many other known examples of
groups with property $(T)$ also have property $F \otimes H$ simply
by verifying the hypotheses of Lemma \ref{lemma:stupid} for those
groups. This is clear for higher rank semisimple algebraic groups
defined over local fields other than $\R$ or $\C$. For lattices in
those groups, it follows from Margulis' superrigidity theorems
\cite{Margulis-Book}. For quotients of lattices in $SP(1,n)$ and
$F_4^{-20}$ by infinite index normal subgroups it follows from
Corlette's superrigidity theorem \cite{Corlette-Annals}.  For any
group where property $(T)$ derives from an action on a simplicial
complex with a spectral condition on links, it can be deduced from
work of Wang, Schoen-Wang or Izeki-Nayatani \cite{Wang, Schoen-Wang,
Izeki-Nayatani}. It seems to us that for most models of random
groups it is at least the case that $(T)$ implies $F \otimes H$ with
probability one. This is most clear in Zuk's model of random group
where the proof of property $(T)$ amounts to proving that the random
group has a properly discontinuous cocompact action on a simplicial
complex with an appropriate spectral condition on links \cite{Zuk}.
For a longer discussion/comparison between models of random groups,
we recommend \cite{Ollivier-Book}.

We believe that non-uniform lattices in semisimple Lie groups with
property $(T)$ will also have property $FH$.  While the case of
$\mathbb Q$-rank $1$ lattices seems quite approachable, the
general case involves some significant difficulties.  We hope to
pursue this elsewhere.

There is a class of groups with property $(T)$ for which we
currently do not know if property $F \otimes H$ holds.  If $G$ is a
semisimple Lie group (as above), there can exist infinite, non-split
central extensions of $G$ by $\Z^d$.  The simplest examples of these
occur when $G$ is the isometry group of an irreducible K\"{a}hler
symmetric space and the universal cover $\tilde G$ is a central
extension of $G$ by $\Z$. Other examples also exist.  It is known
that these groups and their lattices have property $(T)$ by an
argument due to Serre, exposed in \cite{delaHarpe-Valette-Book}.  We
do not know if these groups have property $F \otimes H$. Serre's
argument does not adapt easily to our setting, so it appears that a
new idea is needed for this case.

\section{Reductions to cohomology of cocompact lattices}
\label{section:reductions}

In this section we begin the proof of Theorem \ref{theorem:ft} by
reducing to the following:

\begin{thm}
\label{thm:1-vanishing} Let $G$ be a semisimple Lie group with
property $(T)$ of Kazhdan and no compact factors, and $\G$ a
torsion-free, irreducible, cocompact lattice in $G$. Let
$\pi_1:G{\rightarrow}GL(F)$ be a non-trivial, finite dimensional,
real representation that does not contain the trivial representation
and $\pi_2:\G{\rightarrow}{\cu}(\fh)$ be a (real) unitary
representation on a separable Hilbert space. Then
$H^1(\G,{\pi_1}_{\mid \Gamma}{\otimes}\pi_2,F{\otimes}\fh)=0$.
\end{thm}

\noindent This theorem is proved in section
\S\ref{section:thehardpart} below.

\subsection{First reductions on $G$ and $\Gamma$}
In this section, we show how to reduce to the case of $\G$ an
irreducible torsion-free lattice in a non-compact semisimple group
or $G$ a simple Lie group. The results are not difficult and we
leave the proofs to the reader.

\begin{lem}
\label{lemma:noncompact suffices} Any compact group has property $F
\otimes H$.
\end{lem}

Hence, it suffices to consider the case where $G$ is noncompact.
By a lemma of Selberg, every lattice contains a torsion-free
normal subgroup of finite index \cite{Selberg}. We also have the
following standard lemma.
\begin{lem}\label{lemma:torsion free suffices}
Let $\Gamma'$ be a finite index normal subgroup of $\Gamma$ and
$\rho$ any representation of $\G$ on a real or complex vector
space. Then $H^p(\Gamma', \rho_{\mid_{\Gamma'}})=0$ implies that
$H^p(\Gamma,\rho)=0$.
\end{lem}
\noindent Therefore, we assume from now on that any lattice is
torsion free and pass to finite index subgroups without further
comment.

Recall that a lattice in a semisimple Lie group is irreducible if
it projects densely into all simple factors. If $\G$ is not
irreducible then, possibly after passing to a subgroup of finite
index in $\G$, we have $G=G_1\times G_2$, $\G =\G_1\times \G_2$,
where $\G_i$ is a lattice in $G$ \cite[Definition II.6.5 and
Theorem II.6.7]{Margulis-Book}. So if $\Gamma$ is not irreducible,
after passing to a subgroup of finite index, there is a (split)
short exact sequence of groups $0\To \Gamma_1 \To \Gamma\To
\Gamma_2\To 0$. The resulting long exact sequence in cohomology
contains a portion of the form
\begin{equation}\label{equation:long exact sequence}
\cdots \To H^1(\Gamma_2, \rho) \To H^1(\Gamma, \rho) \To
H^1(\Gamma_1, \rho) \To \cdots.
\end{equation}
Hence, vanishing of $H^1$ for both $\Gamma_1$ and $\Gamma_2$
implies $H^1(\G,\rho)=0$. Induction on the number of irreducible
factors reduces us to the case of irreducible lattices.

When proving Theorem \ref{theorem:ft} for Lie groups, an argument
like the one just given allows us to reduce to the case of $G$ a
simple non-compact Lie group.

\subsection{Decompositions of representations}
\label{subsection:firstfewtricks}

In this subsection we show how reduce to the case $\pi_1$ does not
contain a unitary subrepresentation and when all representations
are real. Recall that a representation which is $\C$ linear is
also $\R$ linear.

\begin{lem}
\label{lemma:realreps} Let $\rho$ be any continuous complex linear
representation of a topological group $D$ on a complex Hilbert
space $W$ and let $\rho|_{\R}$ be the associated real
representation, then $H^1(D,\rho,W)=H^1(D,\rho_{\R},W)$ as real
vector spaces.
\end{lem}

The next observation is a special case of a general decomposition
into irreducibles for representations of the type we consider.  We
do not state the general version here, simply to avoid a
discussion of direct integrals. Again, the lemma stated is not
difficult.

\begin{lem}
\label{lemma:decomp} Let $D$ be a topological group and $\pi$ a
continuous finite dimensional by unitary representation of $D$.
Assume $\pi_1$ is completely reducible and let $({\pi_1},F)=\oplus
(\pi_1^i, F_i)$ be a decomposition into irreducible representations,
then

$$H^1(D, {\pi_1} {\otimes} \pi_2, F{\otimes}\fh)= \sum_i
H^1(D, \pi^i_1 {\otimes} \pi_2, F^i{\otimes}\fh).$$

In particular, if $D$ has property $(T)$, and any $F^i$ is a
unitary $D$ representation, then $\pi_1^i{\otimes}\pi_2$ is
unitary and $H^1(\Gamma, \pi^i_1 {\otimes} \pi_2,
F^i{\otimes}\fh)=0$.
\end{lem}

\noindent{\bf On decomposing $\pi_2$:} It is also possible to
decompose $\pi_2$ into irreducible representations using a direct
integral. In this case, the analogue of Lemma \ref{lemma:decomp}
is not obvious, since solving a cohomological equation in every
irreducible integrand of a direct integral does not necessarily
solve the integral of the cohomological equation. This is already
true for infinite direct sums, see \S\ref{section:gkcohomology}
for more discussion.

\begin{lem}\label{lemma:separable suffices}
Let $D$ be either a countable discrete group or a topological
group containing a countable dense subset.  Then to prove
$H^1(D,\pi)=0$ for all continuous finite dimensional by unitary
representations $\pi$ it suffices to consider
$\pi=\pi_1{\otimes}\pi_2$ where $\pi_2$ is a representation on a
separable Hilbert space.
\end{lem}

\begin{proof}
We prove the lemma first for $D$ a countable discrete group and then
explain the modifications necessary for when $D$ is not discrete.

Let $\pi=\pi_1{\otimes}{\pi_2}$ be a finite dimensional by unitary
representation of $D$ on $F{\otimes}\fh$ for $\fh$ arbitrary and $c$
a cocycle in $H^1(D,\pi_1{\otimes}\pi_2, F{\otimes}\fh)$. We will
show that $c$ takes values in $F{\otimes}\fh_c$ where $\fh_c$ is a
separable Hilbert space  and $F{\otimes}\fh_c$ is $\pi$ invariant.
It then follows that $c$ is a coboundary for $\pi$ if and only if it
is a coboundary for $\pi|_{F{\otimes}\fh_c}$. Consider the set
$D{\cdot}c=\{\pi(d_1)c(d_2)|d_1,d_2{\in}D\}$. This set is countable.
If $\pi$ were a unitary representation, we would take the span of
the closure of $D{\cdot}c$ and consider the restriction of $\pi$ to
this, necessarily separable, subspace. In our setting, a little more
care is required.  Since $F$ is finite dimensional, we can pick a
finite basis, $\{f_1,\ldots, f_k\}$ for $F$.  Every element $\phi$
of $D{\cdot}c$ can now be written as $\phi=\sum_{j=1}^k f_i
{\otimes} \psi^{\phi}_i$ where $\psi^{\phi}_i{\in}\fh$.  It is clear
that $(D{\cdot}c)^{\fh}=\{\psi_i^{\phi} | {\phi}{\in}D{\cdot}c, 1
\leq i \leq k\}$ is a countable subset of $\fh$.  Let $\fh_c$ be the
closure of the span of $(D{\cdot}c)^{\fh}$, then  it is clear that
$\fh_c$ is a closed separable Hilbert subspace of $\fh$, that
$\fh_c$ is $D$ invariant and that $c$ takes values in
$F{\otimes}\fh_c$. If $c$ is a coboundary, then $c(d)=\pi(d)v-v$. If
such a $v$ exists, it clearly can be chosen in $F{\otimes}\fh_c$.
Therefore $c$ is a coboundary as a cocycle over $\pi$ if and only if
$c$ is a coboundary as a cocycle over $\pi|_{F\otimes\fh_c}$.

If $D$ is a topological group with countable dense subset $D_0$,
we need only replace $D{\cdot}c$ by the countable set
$D_0{\cdot}c=\{\pi(d_1)c(d_2)|d_1,d_2{\in}D_0\}$ and the same
argument works. Continuity of $\pi$ and $c$ imply that, in this
case, the closure of the span $D_0{\cdot}c$ equals the closure of
the span of $D{\cdot}c$.  This then implies that $F{\otimes}\fh_c$
constructed as above is a $D$ invariant closed, separable subspace
of $F{\otimes}\fh$.
\end{proof}

\subsection{Reductions via superrigidity theorems}

In this subsection we recall a form of the
Margulis-Corlette-Gromov-Schoen superrigidity theorems. Given a
lattice $\G$ in $G$, and a representation
$\pi:\Gamma{\rightarrow}GL(F)$, we say that $\pi$ {\em almost
extends} to a continuous representation of $G$ if there exist
representations $\pi_1:G{\rightarrow}GL(F)$ and
$\pi_2:\Gamma{\rightarrow}GL(F)$, where $\pi_2$ has bounded image
and the images of $\pi_1$ and $\pi_2$ commute so that
$\pi(\g)=\pi_1(\g)\pi_2(\g)$.

\begin{thm}
\label{theorem:superrigidity} Let $G$ be a semisimple Lie group with
no compact factors and property $(T)$ of Kazhdan.  Let $\Gamma$ in
$G$ be a lattice and $\pi:\G{\rightarrow}GL(F)$ a finite dimensional
representation. Then $\pi$ almost extends to a continuous
representation of $G$.
\end{thm}


\begin{proof}
This is a consequence of the strongest form of the superrigidity
theorems. All the ingredients needed are assembled by Starkov in
\cite{Starkov}, though the theorem is not stated in quite this
generality there.

For $G$ of higher rank and $\G$ irreducible, this is proved in
\cite{Margulis-Book}, though it is not stated explicitly there. As
noted in \cite{Fisher-Margulis0}, it follows easily from from
Lemma VII.5.1 and Theorems VII.5.15 and VII.6.16 of
\cite{Margulis-Book}. However, when $G$ contains rank one factors
and $\G$ (or some finite index subgroup) projects to a lattice in
a rank one factor one needs to use Corlette's work
\cite{Corlette-Annals} to prove the analogue of \cite[Theorem
VII.5.15]{Margulis-Book}.  The proof of Theorem VII.6.16 is even
more involved, as one needs to first prove arithmeticity using
work of Corlette and Gromov-Schoen and then deduce the theorem
from cohomology vanishing theorems for finite dimensional
representations of $\G$.  As noted in \cite{Starkov} this follows
by assembling results from \cite{Raghunathan-1-vanishing} and
\cite{Raghunathn-1-vanishing-nonuniform}.  At various points in
the argument, one needs to pass to finite index subgroups, but
this is not required in the statement of results by use of
\cite[Lemma VII.5.1]{Margulis-Book}.
\end{proof}

\begin{cor}
\label{corollary:extendabletensorunitary} Let $G$ and $\Gamma$ be as
in Theorem \ref{theorem:superrigidity}.  Let
$\pi=\pi_1{\otimes}\pi_2$ be any $\Gamma$ representation as in
Definition \ref{definition:ft}.  Then we can assume $\pi_1$ extends
to $G$.
\end{cor}

\begin{proof}
It is easy to check that if a representation $\pi$ of $\G$ can be
written as $\pi(\g)=\pi'(\g)\pi''(\g)$ where $\pi'$ and $\pi''$
commute, then $\pi$ can be realized as a tensor product of
representations $\pi'_0{\otimes}\pi''_0$. So, if $\pi_1$ does not
extend to $G$, then as a consequence of Theorem
\ref{theorem:superrigidity}, we can write
$\pi_1=\pi_1'{\otimes}\pi_2'$ where $\pi_1'$ extends and $\pi_2'$
has bounded image. Regrouping $\pi$ as
$\pi_1'{\otimes}(\pi_2'{\otimes}\pi_2)$ proves the corollary.
\end{proof}

\subsection{Reduction to cocompact lattices}

The purpose of this subsection is to show how one reduces the case
of connected groups to the case of cocompact lattices. In the
following lemma, we call a lattice in a product $G_1{\times}G_2$
\emph{weakly irreducible} if both projections on factors are dense.

\begin{lem}
\label{lemma:restriction} Let $G$ be a locally compact group and
assume $G{\times}G$ admits a cocompact weakly irreducible lattice
$\Gamma$. Then if $\Gamma$ has property $F \otimes H$ so does $G$.
\end{lem}

\begin{proof}
Let $\G$ be a weakly irreducible lattice in $G{\times}G$ and
$(\pi,F{\otimes}\fh)$ a finite dimensional by unitary representation
of $G$. If $p_1$ is the projection of $\G$ on the first factor of
$G{\times}G$, then $\pi{\circ}p_1$ is also clearly a finite
dimensional by unitary $\G$ representation. We will show that
$H^1(\G, \pi{\circ}p_1,F{\otimes}\fh)=0$ implies $H^1(G,\pi,
F{\otimes}\fh)=0$. Any continuous cocycle
$c:G{\rightarrow}F{\otimes}\fh$ restricts to a $p_1(\Gamma)$
cocycle.  A coboundary for $\Gamma$ is a map of the form
$b(\g)=\pi{\circ}p_1(\g)v-v$. Note that continuity of $\pi$ and the
formula for $b$ immediately imply that $b$ is continuous in the
induced topology on $p_1(\G)<G$ and so $b$ extends continuosly to a
$G$ coboundary $b(g)=\pi(g)v-v$. If $H^1(\G,
\pi{\circ}p_1,F{\otimes}\fh)=0$, then $c{\circ}p_1(\g)=b(\g)$. Since
both sides extend continuously to functions on $G$, they are equal
on $G$ as well, so $c(g)=b(g)$ and $H^1(G,\pi,F{\otimes}\fh)=0$ as
desired.

It is also possible to write this proof in terms of fixed points
for affine actions.  From this point of view, the idea is that for
a continuous affine $G$ action, a fixed point for a dense subgroup
is necessarily a fixed point for $G$.
\end{proof}

{\noindent}{\bf On Lemma \ref{lemma:restriction}:} \begin{enumerate}

\item The lemma is a fact about $H^1$ and has no obvious analogue
for higher degree cohomology. The key point is that in a
continuous $G$ representation, a $1$-coboundaries for a dense
subgroup extends continuously to a $1$-coboundary for $G$. This is
not clear, and seems unlikely, for $k$-coboundaries when $k>1$.

\item We recall that for every group $G$ as in Theorem
\ref{theorem:ft}, there is a weakly irreducible cocompact lattice in
$G{\times}G$.  This follows from our definition of irreducible
lattice (note that we do not require that the projection to other
quotients of $G{\times}G$ be dense) and Borel's construction of
lattices in \cite{Borel}.  To construct the lattice, one can simply
construct an irreducible lattice (in the standard sense) in $G_i
\times G_i$ for each simple factor $G_i$ of $G$. Due to prior
reductions, for our argument we need only construct irreducible
cocompact lattices in $G{\times}G$ where $G$ is simple.
\end{enumerate}

\subsection{Proof that Theorem \ref{thm:1-vanishing} implies
Theorem \ref{theorem:ft}}

\begin{proof}
By Lemma \ref{lemma:noncompact suffices}, Lemma \ref{lemma:torsion
free suffices} and the discussion around Equation \ref{equation:long
exact sequence}, we see that we may assume that $G$ has no compact
factors and that $\G$ is irreducible and torsion free. By Lemma
\ref{lemma:realreps} it suffices to consider the case where
$\pi=\pi_1{\otimes}\pi_2$ where $\pi_1$ is a real linear
representation and $\pi_2$ is a unitary representation over $\R$. By
Lemma \ref{lemma:separable suffices}, it suffices to consider the
case where $\fh$ is a separable Hilbert space.

By Theorem \ref{theorem:superrigidity} it suffices to prove
Theorem \ref{theorem:ft} for representations which are of the form
$\pi_1{\otimes}\pi_2$ where $\pi_1:G{\rightarrow}GL(F)$ is a real
linear representation and $\pi_2$ is a unitary representation of
$G$ or $\Gamma$. By Lemma \ref{lemma:decomp} and the fact that $G$
representations are completely reducible, it suffices to consider
the case where $\pi_1$ does not contain the trivial
representation.

Lemma \ref{lemma:restriction} and the second remark following that
lemma show how Theorem \ref{theorem:ft} for cocompact lattices
implies Theorem \ref{theorem:ft} for connected groups.
\end{proof}

\section{Vanishing theorems for uniform lattices}
\label{section:thehardpart}
 In this section we prove Theorem
\ref{thm:1-vanishing}. Actually, we prove a more general result than
Theorem \ref{thm:1-vanishing} which includes some cases where $G$
has factors of the form $\mathrm{SO}(1,n)$ or $\mathrm{SU}(1,n)$,
see \S\ref{subsection:sufficienttovanish}. For this reason, until
specified further, in this section $G$ is a semisimple Lie group
with no compact factors and $\G \subset G$ is a cocompact lattice.
We also derive some results on higher cohomology groups in the case
that the finite dimensional representation is complex instead of
real. See \S \ref{subsection:higherdegrees} for specific statements.
For this reason, much of the development in this section treats
$p$th cohomology groups rather than just $H^1$.

The main technical tool is an estimate on the Laplacian on smooth
differential forms with values in a vector bundle. Assume that
$\Gamma$ is torsion free and let $K$ be a maximal compact subgroup
of $G$. We build the locally symmetric Riemannian manifold $M =
\Gamma \backslash G/K$ and the locally constant vector bundle
$\mathcal{E}(\pi) \To M$ associated to the representation $\pi$
and the principal $\Gamma$-bundle $\Gamma\To G/K \To M$. Note that
the fiber of this bundle is the infinite dimensional vector space
$E=F\otimes\mathcal{H}$.

We construct an inner product in the fibers of $\mathcal{E}(\pi)$,
and therefore inner products $(\cdot ,\cdot )$ on $\bigwedge^p
T^*M\otimes \mathcal{E}(\pi)$, the spaces of smooth
$\mathcal{E}(\pi)$-valued differential $p$-forms. Using these inner
products, we define an adjoint operator $\delta$ to the exterior
derivative operation $d$ and hence a Laplace operator $\triangle = d
\delta + \delta d$. The estimate we obtain is the following.
\begin{thm}\label{estimate on forms} Let $G, \G$ and $\pi$ be as
in \ref{thm:1-vanishing}. There is a positive constant
$C=C(G,\pi_1)$ such that for any smooth differential $1$-form,
$\eta$, with values in $\mathcal{E}(\pi)$  we have
$$( \triangle \eta, \eta )\geq C(\eta, \eta).$$
\end{thm}

\noindent We will need this result itself for our geometric
cocycle superrigidity theorem. It will be translated into an
energy estimate and will play a crucial role by controlling a type
of heat flow on a special set of mappings
\cite{Fisher-Hitchman2,Fisher-Hitchman1}.

We now outline the rest of the section. In
\S\ref{subsection:cohomology change} we discuss some cohomology
isomorphisms which reduce the problem to considering de Rham
cohomology of smooth $p$-forms with values in a vector bundle. In
\S\ref{subsection:laplacians} we introduce the necessary operations
on forms to define the Laplacian on a vector bundle. In
\S\ref{subsection:sufficienttovanish} we discuss how vanishing
cohomology can be deduced from an estimate on this Laplacian. The
material in these three sections is mostly standard, though perhaps
not well-known in our setting.

In \S\ref{Bochner formula section}, we derive a useful Bochner-type
formula for the Laplacian on our particular vector bundle. While our
final estimate is similar to that of Matsushima and Murakami, the
derivation follows a different conceptual outline. We then use this
formula and some Lie-theoretic computations of Raghunathan to obtain
the required estimates and deduce vanishing results. We discuss
first cohomology in \S \ref{subsection:endgame} and higher degrees
in \S\ref{subsection:higherdegrees}. In \S
\ref{subsection:flatlaplacian}, we give a direct computation of the
flat Laplacian.

\subsection{Some Cohomology
Isomorphisms}\label{subsection:cohomology change}

Let $K$ be a maximal compact subgroup of $G$.  Let $M = \Gamma
\backslash G/K$, a compact locally symmetric manifold, and form
the locally constant vector bundle $\mathcal{E}(\pi) \To M$
associated to the representation $\pi$ and the principal
$\Gamma$-bundle $\Gamma\To G/K \To M$. Note that the canonical
fiber of this bundle is $E=F\otimes\mathcal{H}$. Let
$\mathcal{E}(\pi)_{loc}$ denote the system of local coefficients
on $M$ defined by the bundle $\mathcal{E}(\pi)\To M$.

\begin{prop} \label{cohom iso 1} $H^*(\Gamma, \pi) = H^*(M, \mathcal{E}(\pi)_{loc} )$
\end{prop}

\noindent This isomorphism holds on the level of cochains
\cite{Brown-cohomofgps}. Let $\check{H}^*(M, \mathcal{E})$ denote
the \v{C}ech cohomology of $M$ with coefficients in the sheaf
$\mathcal{E}$ of locally constant sections of $\mathcal{E}(\pi)$.

\begin{prop}\label{cohom iso 2} $H^*(M, \mathcal{E}(\pi)_{loc} )=  \check{H}^*(M, \mathcal{E})$
\end{prop}

\noindent One can see this by inspecting the definitions. For more
on sheaf cohomology, see for example \cite{Warner-book}. Let
$H^*_{\text{de Rham}}(M, \mathcal{E}(\pi))$ denote the de Rham
cohomology of smooth $\mathcal{E}(\pi)$-valued differential forms on
$M$. The exterior derivative $d$ on $\mathcal{E}(\pi)$-valued
$p$-forms is defined as follows. For a form of the type $\omega =
\eta\otimes e$, where $\eta$ is an ordinary $p$-form on $M$ and $e$
is a section of $\mathcal{E}(\pi)$, we set
$$d\omega = d\eta \otimes e,$$
and extend linearly. Then one can define de Rham cohomology in the
usual fashion.

\begin{prop}\label{cohom iso 3} $\check{H}^*(M, \mathcal{E}) \cong H^*_{\text{de Rham}}(M,
\mathcal{E}(\pi))$.\end{prop}

\noindent This is discussed by Mok \cite{Mok}. The key point is
that a version of the Poincar\'{e} lemma holds. Using these
isomorphisms, we restrict ourselves to the study of de~Rham
cohomology.

\subsection{The Laplacian on a Euclidean Vector
Bundle}\label{subsection:laplacians}

In this section we recall some facts about Riemannian vector
bundles. First, we recall the definitions of some standard
operations which lead to a definition of the inner product and the
Laplacian on the space of bundle-valued $p$-forms. Then we recall
the Bochner-Weitzenb\"{o}ck formula for bundle-valued differential
forms. A reference for this material is \cite{Eells-Lemaire}.

Throughout this section, we consider a vector bundle $V$ over a
compact, oriented, Riemannian manifold $N$. We assume also that this
bundle is Riemannian, that is, $V$ is equipped with a connection
$\nabla$ and a fiberwise Euclidean structure $(\ ,)$ which are
compatible in the sense that for any two sections $a,b$ of $V$ and
any continuous vector field $X$ on $N$,
$$X\cdot (a,b) = (\nabla_X a, b) + (a, \nabla_X b).$$
We denote by $A^p(V)$ the space of smooth $p$-forms on $N$ with
values in $V$, i.e. the space of smooth sections of the bundle
$\bigwedge^p T^*N\otimes V$ over $N$. Given the pair $V,\nabla$,
one can define the exterior derivative operator $d$ for $\sigma
\in A^p(V)$ by
\begin{multline*}
    d\sigma(X_1,\ldots,X_{p+1}) = \sum_{i=1}^{p+1}
    (-1)^{i+1}\nabla_{X_i}(\sigma(X_1,\ldots,\hat{X_i},\ldots,X_{p+1}))\\
    +\sum_{i<j} (-1)^{i+j} \sigma([X_i,X_j],X_1,\ldots, \hat{X_i}, \ldots,
    \hat{X_j}, \ldots, X_{p+1}),
\end{multline*}
where the hat means that that term is omitted. Note that $d^2$ is
not necessarily zero. In fact, $d^2=0$ if and only if the
connection $\nabla$ is flat.

In the case of $\mathcal{E}(\pi)$, the differential $d$ which
computes de~Rham cohomology is constructed by this process,
$d^2=0$ and the bundle is flat.

\subsubsection{The inner product on $V$-valued forms} The Hodge
star operation for ordinary $p$-forms on $N$ induces a similar
operation on $V$-valued forms which we shall also denote by $*$.
For $\omega =\sum \eta_I \otimes e_I$ where $I$ is a multi-index
of length $p$, $\eta_I$ is an ordinary $p$-form and $e_I$ is a
section of $V$, we have
\begin{equation*}
    * \omega = \sum (*\eta_I) \otimes e_I.
\end{equation*}

Let $V^*$ be the dual vector bundle to $V$. The Euclidean structure
in the fibers of $V$ defines an isomorphism $\#: V \To V^*$ which is
determined by the condition that for $u_x, v_x$ in the fiber of $V$
over $x\in N$
$$\langle u_x, \# v_x\rangle= (u_x,v_x)_x,$$
where $\langle \cdot, \cdot\rangle$ denotes the natural pairing of
$V$ and $V^*$. We note that the operations $*$ and $\#$ commute.

We define a wedge product for $V$-valued forms. Let $\theta=
\sum_{|I|=p} \theta_I\otimes e_I$ be an $V$-valued $p$-form and
$\eta=\sum_{|J|=q} \eta_J\otimes f_J$ an $V^*$-valued $q$-form. We
define an ordinary $(p+q)$-form, $\theta^t\wedge\eta$, on $N$ as
follows
\begin{equation*}
    \theta^t\wedge \eta = \sum_I\sum_J
    \langle e_I,f_J\rangle\theta_I\wedge \eta_J.
\end{equation*}
The transpose notation is used since the forms live in different
bundles and this is not the standard wedge product.

We now define the inner product of a pair of $V$-valued $p$-forms.
For two such forms $\omega_1, \omega_2$, we set
\begin{equation}\label{inner product definition}
    (\omega_1, \omega_2) = \int_N \omega_1^t\wedge *\#\omega_2 =
     \int_N \omega_1^t\wedge \# *\omega_2.
\end{equation}
We use the same symbol here for the inner product of forms, the
inner product in the fibers of $V$ and the induced inner products
in various associated bundles. Which structure we mean will be
clear from context.

The integrand can be evaluated in local coordinates. Pick a point $x
\in N$. In a sufficiently small neighborhood $U$ of $x$, we can
choose a basis of one forms $\theta^1, \ldots, \theta^n$ on $N$
which are orthonormal at each point of $U$ and positively oriented.
Write $V$-valued $p$-forms $\omega_1,\omega_2$, in coordinates
\begin{eqnarray*}
  \omega_1 = & \sum_{i_1<\ldots<i_p} \theta^{i_1}\wedge \ldots \wedge
  \theta^{i_p} \otimes a_{i_1,\ldots,i_p} \\
    \omega_2 = & \sum_{i_1<\ldots<i_p} \theta^{i_1}\wedge \ldots \wedge
  \theta^{i_p} \otimes b_{i_1,\ldots,i_p} \\
\end{eqnarray*}
where $a_{i_1,\ldots,i_p}$ and $b_{i_1,\ldots,i_p}$ are sections
of $V$. Then at all points $z$ of $U$,
\begin{equation}\label{inner product integrand locally}
 \omega_1^t\wedge
*\#\omega_2(z) = \sum_{i_1<\ldots<i_p}( a_{i_1,\ldots,i_p},
b_{i_1,\ldots,i_p})(z).
\end{equation}

\subsubsection{The operators $\delta$ and $\triangle$} We now define the adjoint of
exterior differentiation and the Laplacian. For any $V$-valued
$p$-form $\theta$ on $N$, we define
$$\delta\theta = (-1)^p *^{-1} \#^{-1} d \# * \theta. $$
\begin{prop} The operator $\delta$ is the formal adjoint of $d$ with respect to
the inner product $(\cdot, \cdot)$ on $V$-valued $p$-forms. That
is, for two $V$-valued $p$-forms $\theta$ and $\omega$:
\begin{equation*}
    (\delta \theta, \omega) = (\theta, d\omega).
\end{equation*}
\end{prop}

For a smooth $V$-valued $p$-form on $N$ we define the Laplacian
by:
\begin{equation*}
    \triangle \theta = d\delta \theta + \delta d\theta.
\end{equation*}
\noindent This is a self-adjoint elliptic differential operator of
second order on (sections of) $V$. We note that it is also a
non-negative operator, since
    $(\triangle \theta, \theta) = (d\theta, d\theta) + (\delta\theta,
    \delta\theta) \geq 0$.
\noindent The formal adjoint $\nabla^*$ of $\nabla$ is defined by
the property that $(\nabla \eta, \theta) = (\eta,
\nabla^*\theta)$.

\subsubsection{The traditional Bochner-Weitzenb\"{o}ck formula}
Here we recall the Bochner-Weitzenb\"{o}ck formula for $V$-valued
$p$ forms. For a development from first principles, see \cite[pages
3-13]{Eells-Lemaire}.

Let $R^V$ denote the curvature tensor of the connection $\nabla$ on
$V$ and $R$ the Riemannian curvature tensor of $N$. We define the
Ricci operator, $S$, on $V$-valued $p$-forms as follows. Let $X_1,
\ldots, X_p$ be continuous vector fields on $N$, $x$ a point in $N$,
and $e_1, \ldots, e_n$ an orthonormal basis of $T_xN$. Then for a
$V$-valued $p$-form $\sigma$, $S(\sigma)$ is the $V$-valued $p$-form
given at $x$ by
\begin{multline}\label{ricci operator}
    S(\sigma)_x(X_1, \ldots, X_p) = \\
    = \left\{ \begin{array}{cc}
        0, & \text{if $p=0$}, \\
        \sum_{k=1}^p\sum_{s=1}^n (-1)^k
        (R^V(e_s,X_k)\sigma)(e_s, X_1, \ldots, \widehat{X_k}, \ldots,
        X_p), & \text{if $p\geq 1$}. \\
    \end{array}\right.
\end{multline}

\begin{thm}[Eells-Lemaire, p. 13]
Let $\sigma$ be a $V$-valued $p$-form on $N$. Then at each point $x$
of $N$ we have
\begin{equation}\label{equation:pointwisebochnerformula}
    \dfrac{1}{2} \triangle (\sigma, \sigma)(x)= (
    \triangle \sigma, \sigma)(x) - ( \nabla\sigma,
    \nabla\sigma)(x) - ( S(\sigma), \sigma)(x).
\end{equation}
\end{thm}

\noindent This is a pointwise formula, the inner products are the
Euclidean product in the fiber over $x$ taken at the values of the
forms at $x$. Integrating over $N$ yields:

\begin{cor}\label{integrated bw formula}
\begin{equation}
\label{equation:integratedbw}
    \int_N ( \triangle \sigma, \sigma) dvol_N = \int_N
    (\nabla^*\nabla \sigma, \sigma) dvol_N + \int_N (
    S(\sigma), \sigma) dvol_N.
\end{equation}
\end{cor}

\noindent The corollary follows as the left hand side of equation
\ref{equation:pointwisebochnerformula} integrates to zero by the
divergence theorem, and $\nabla^*$ is the formal adjoint  of
$\nabla$.

\subsection{A Sufficient Condition for Vanishing}
\label{subsection:sufficienttovanish}

In this section we show that an estimate of the type in Theorem
\ref{estimate on forms} implies vanishing of cohomology. The
argument is close to one due to Mok, compare to the proof of
\cite[Proposition 1.3.1]{Mok}.
\begin{prop}\label{sufficient condition for vanishing}
Suppose that there exists a positive constant $C$ such that for
all smooth $\mathcal{E}(\pi)$-valued $p$-forms $\theta$ on $M$
\begin{equation}
\label{Laplaceestimate}
    (\triangle\theta, \theta) \geq C (\theta, \theta).
\end{equation}
Then $H^p_{\text{de Rham}}(M, \mathcal{E}(\pi)) =0$.
\end{prop}

\begin{proof} In order to give the proof, we need a more general
viewpoint on our forms to employ some functional analysis. Let
$A^p(M,\mathcal{E}(\pi))$ be the Hilbert space of
$\mathcal{E}(\pi)$-valued differential forms $\eta$ on $M$ for which
$\eta, d\eta$ and $\delta\eta$ satisfy an $L^2$ condition with
respect to our inner product $(\cdot, \cdot)$. Smooth $p$-forms are
a dense in $A^p(M,\mathcal{E}(\pi))$. By a simple approximation
argument, condition \ref{Laplaceestimate} of Proposition
\ref{sufficient condition for vanishing} holds for all forms in
$A^p(M,\mathcal{E}(\pi))$ if it holds for all smooth forms. We need
the following lemma from functional analysis

\begin{lem}[H\"{o}rmander \cite{Hormander-book}, Lemma 4.1.1, page 78]
\label{theorem:hormander} Let $\phi:H_1\To H_2$ be a densely
defined linear operator between two Hilbert spaces with adjoint
$\phi^*$. Suppose that $H_3\subset H_2$ is a closed linear
subspace which contains the range of $\phi$. Then the range of
$\phi$ is dense in $H_3$ if and only if for some constant $C>0$
\begin{equation}
\label{equation:hormander}
 ||v||_{H_2} \leq C ||\phi^*v||_{H_1}
 \end{equation}
for all $v$ in the intersection of $H_3$ with the domain of
$\phi^*$.
\end{lem}

We apply this result to the the exterior derivative map
$$d:A^{p-1}(\mathcal{E}(\pi)) \To \ker d \subset
A^p(\mathcal{E}(\pi)).$$ Note that for an element $\eta$ of $\ker
d$ we have $(\triangle \eta,\eta) = (\delta \eta, \delta \eta)$,
so equation \ref{equation:hormander} follows from equation
\ref{Laplaceestimate}.

Since any smooth closed $\mathcal{E}(\pi)$-valued $p$-form $\eta$
lies in the kernel of $d$, we see that there is a solution to the
equation $du =\eta$. It remains to be show that $u$ can be chosen
smooth. This is true for the solution having minimal $L^2$-norm.

To see that a solution of minimal norm exists, let
$t=\inf\{||u||_2 \mid du = \eta\}$ and consider a sequence $u_i$
of elements of $A^p(M,\mathcal{E}(\pi))$ with $du_i$ approaching
$t$. The ball of $L^2$-forms having norm bounded by $2t$ is weak-*
compact, so we may extract a weak-* subsequential limit of the
$u_i$. By lower semicontinuity of the $L^2$ norm, we see that $$
||u||_2 \leq \liminf ||u_i||_2 = t. $$ We claim that this limit
$u$ is also a solution to $du=\eta$. Let $\psi$ be a smooth test
form. We have $(\eta, \psi) = (du_i, \psi ) = (u_i, \delta \psi)
\rightarrow (u, \delta \psi) = (du, \psi)$. Hence $du=\eta$.

We now show that $\delta u = 0$. Suppose that $\delta u =g$ is not
zero. Then there exists an $h$ such that $(g,h) = (\delta u, h) =
(u, dh) =: b > 0.$ Notice this is also true when we replace $h$ by
$r h$ where $r$ is a positive real number. Let $a=(dh,dh)>0$ and
pick a positive real number $r < 2a/b$. Let $w= u - r dh$. Then
$dw = du = f$, and $(w,w) = (u,u) - 2r(u,dh) + r^2(dh,dh) =(u,u) -
r(2a-br) <(u,u)$. This contradicts the minimality of the norm of
$u$, so we conclude that $\delta u= 0$.

Now note that the form $\nu=\delta\eta$ is a smooth form, and by
construction $u$ satisfies
\begin{equation*}
    \triangle u = (d\delta + \delta d) u = \delta\eta =
    \nu.
\end{equation*}
So the elliptic regularity theorem implies that $u$ is a smooth
form. Thus we see that every smooth closed $\mathcal{E}(\pi)$-valued
$p$-form bounds.
\end{proof}

\subsection{A variant of the Matsushima-Murakami Bochner
formula}\label{Bochner formula section}

Our goal in this section is to develop a variant of the
Bochner-Weitzenb\"{o}ck formula for $\mathcal{E}(\pi)$-valued
$p$-forms on $M$ which is analogous to the one found by Matsushima
and Murakami \cite{Matsushima-Murakami}. As with all formulae of
this type, this formula will be a computation of the difference
between two second order operators on $\mathcal{E}(\pi)$-valued
$p$-forms and the difference is a zeroth order (algebraic) operator
related to curvature. Then we analyze this algebraic operator to
prove Theorem \ref{estimate on forms}.

One could give a development directly mirroring that of
\cite{Matsushima-Murakami}. Instead, we give an alternate derivation
which we feel is more intuitive and makes the relation to other
Bochner-type formulae more explicit.

The key observation of Matsushima and Murakami in
\cite{Matsushima-Murakami} is that the relevant vector bundle can
be constructed in two distinct ways, that is, there are two
(isomorphic) natural vector bundles associated to the given data.
This allows us to compare the natural differential operators
arising from the two constructions.

In \S\ref{subsection:bundleisos}, we exhibit two distinct vector
bundle structures on $\mathcal{E}(\pi)$. In
\S\ref{subsection:metric}, we define a preferred Riemannian metric
on $E(\pi)$. We discuss the computation of the connections
compatible with the metric in the two bundle structures in
\S\ref{connections}, and then, in
\S\ref{subsection:explicitformula}, we derive our version of the
Matsushima-Murakami Bochner formula in terms of the representation
$\pi$ and the structure of $G$ (or, rather, its Lie algebra
$\mathfrak{g}$). We deduce the required estimates and the vanishing
theorems which follow in \S\ref{subsection:endgame} and
\S\ref{subsection:higherdegrees}.

\subsubsection{Some bundle isomorphisms}
\label{subsection:bundleisos}

Recall that $\Gamma < G$ is a torsion free lattice and $K < G$ is
a maximal compact subgroup of $G$. We are working with the
situation where $K$ acts on $G$ by right multiplication, $\Gamma$
acts on $G$ by left multiplication and our compact manifold is $M
= \Gamma \backslash G/K$. In fact, $G$ is a right principal
$(\Gamma \times K)$-bundle over $M$, where $\Gamma \times K$ acts
by the rule \begin{equation} \label{equation:rightaction} g\cdot
(\gamma, k) = \gamma^{-1}gk.\end{equation} We denote the quotient
map by $q:G\To M$.

\noindent{\bf On Equation \ref{equation:rightaction}:} Equation
\ref{equation:rightaction} may surprise readers not familiar with
standard conventions on principal bundles. The action described is a
right action, or equivalently, an action of the opposite group. This
convention is used in the literature on connections on principal
bundles and we use it here to remain consistent with the
\cite{KobayashiNomizuI}. This convention is not typically used in
the theory of flat bundles. However, since the opposite group is
canonically isomorphic to the original group, proving the vanishing
result for the opposite group of $\G$ is equivalent to the same
result for $\G$.\\

We are interested in two representations $\sigma$ and $\sigma'$ of
$\Gamma\times K$ on $E=F\otimes \mathcal{H}$, defined as follows.
\begin{eqnarray*}
    \sigma(\gamma,k) v = &(\pi_1(\gamma)\otimes \pi_2(\gamma))v,\\
    \sigma'(\gamma,k) v = &(\pi_1(k)\otimes \pi_2(\gamma))v.
\end{eqnarray*}
We build the vector bundles over $M$ associated to the right
principal bundle $G\To M$ and these representations, respectively
$E(\sigma)$ and $E(\sigma')$. Recall that $E(\sigma)$ is the
quotient of $G\times E$ by the right $\Gamma\times K$-action
$$ (g,v)\cdot (\gamma,k) = (\gamma^{-1}gk, \sigma(\gamma,k)^{-1}v)
= (\gamma^{-1}gk, (\pi_1(\gamma)^{-1}\otimes \pi_2(\gamma)^{-1}) v
),$$ and $E(\sigma')$ is the quotient of $G\times E$ by the
$(\Gamma\times K)$-action
$$ (g,v) \cdot (\gamma, k) =(\gamma^{-1}gk, \sigma'(\gamma,k)^{-1}v)
= (\gamma^{-1}gk, (\pi_1(k)^{-1}\otimes \pi_2(\gamma)^{-1}) v ).$$

\noindent We shall denote the quotient maps here by $Q:G\times
E\To E(\sigma)$ and $Q':G\times E \To E(\sigma')$.  Since
$\sigma(\g,k)=\pi(\g)$ for all $(k,\g) \in K{\times}\G$, we have:
\begin{lem}\label{lifted identification of bundles}
  As vector bundles over $M$, $\mathcal{E}(\pi)\cong E(\sigma)$.
\end{lem}

The following observation is a more general version of
\cite[Proposition 3.1]{Matsushima-Murakami}.

\begin{lem}
  As vector bundles over $M$, $E(\sigma)\cong E(\sigma')$.
\end{lem}

\begin{proof} Consider the mapping $A:G\times E\To G\times E$ defined
by
$$A(g,v) = (g, (\pi_1(g)^{-1}\otimes\mathrm{Id}_{\mathcal{H}})v).$$
This is a bijective mapping which intertwines the $\G{\times}K$
actions inducing the quotient mappings $Q$ and $Q'$. For each $g \in
G$ the map $v\rightarrow A(g,v)$ is a linear isomorphism $E\To E$,
so $A$ induces a vector bundle isomorphism $A':E(\sigma)\To
E(\sigma')$.\end{proof}

\noindent We will use these isomorphisms as identifications
without further comment.

\subsubsection{The Euclidean structure on $\mathcal{E}(\pi)$}
\label{subsection:metric}

We now construct a Euclidean structure on $\mathcal{E}(\pi)$. We
begin with a second important observation of Matsushima and Murakami
in \cite{Matsushima-Murakami}. Let $\mathfrak{g}$ be the Lie algebra
of $G$ with a Cartan decomposition $\mathfrak{g} =
\mathfrak{k}\oplus\mathfrak{p}$ relative to $K$. We abuse notation
and use $\pi_1$ to denote both the representation of $G$ on $F$ and
the induced representation of $\mathfrak{g}$ on $F$.
\begin{prop}[\cite{Matsushima-Murakami}, Lemma 3.1]\label{F-admissible metric}
There is an inner product $(\cdot, \cdot )_F$ on $F$ which satisfies
\begin{itemize}
    \item $( \pi_1(X)u, v)_F = -(u,\pi_1(X)v)_F$ for $X \in
    \mathfrak{k}$, and
    \item $( \pi_1(X)u, v)_F = (u,\pi_1(X)v)_F$ for $X \in
    \mathfrak{p}$.
\end{itemize}
\end{prop}
\noindent The first condition means that $(\cdot , \cdot)_F$ is
invariant under $\pi_1(K)$.

Since $\pi_2$ is a unitary representation of $\Gamma$ on
$\mathcal{H}$, each $\pi_2(\gamma)$ preserves the Hilbert product
$(\ , )_{\mathcal{H}}$ on $\mathcal{H}$. We form an inner product on
$E=F\otimes\mathcal{H}$ as the tensor product of the above two
structures. That is, let $(a\otimes c, b\otimes d)_E =
(a,b)_F\cdot(c,d)_{\mathcal{H}}$ and extend linearly. By the
properties in Proposition \ref{F-admissible metric}, this inner
product is clearly invariant by the action inducing the quotient
$Q':G\times E\To E(\sigma')$, and so induces a Euclidean structure
on $E(\sigma') = \mathcal{E}(\pi)$ which we denote hereafter by $(\
, )$ and call the \emph{canonical} metric on $\mathcal{E}(\pi)$.

\subsubsection{Two linear connections}\label{connections}

The principal bundle $G\To M$ carries a natural connection. The
subgroup $K<G$ induces a Cartan decomposition
$\mathfrak{g}=\mathfrak{k}\oplus\mathfrak{p}$ of the Lie algebra
$\mathfrak{g}$ of $G$, where $\mathfrak{k}$ is the Lie subalgebra
corresponding to $K$. We make a choice of horizontal subspace $H_g
\subset T_gG$ by $H_g = dL_g(\mathfrak{p})$, where $L_g$ denotes
left multiplication by $g$. This family $\{ H_g \}_{g\in G}$
defines a connection in $G\To M$.

This connection in $G\To M$ induces connections $\nabla$ and
$\nabla'$ in the vector bundles $E(\sigma)$ and $E(\sigma')$. Note
that we have associated representations, bundles and connections
for all types of tensors over $M$ with values in one of these two
bundles by the same methods.

Recall the computation of such an induced connection
(\cite{KobayashiNomizuI} page 116). We state the result for
$E(\sigma)$ but the analogous procedure works for $E(\sigma')$.
There is a bijective correspondence between sections $f: M \To
E(\sigma)$ and $E$-valued functions $f^*:G\To E$ such that
\begin{equation}\label{equivariance}
  f^*(g\cdot (\gamma,k)) = \sigma(\gamma,k) f^*(g).
\end{equation}
This correspondence is given in one direction by lifting $f$ to a
section of the trivial bundle $G\times E \To G$ and taking the
second coordinate. Conversely, a function satisfying equation
\ref{equivariance} defines the second coordinate of a section
$G\To G\times E$ which descends via the quotient map $Q: G\times E
\To E(\sigma)$ since it is $\G{\times}K$ invariant.

Now we have the following result.
\begin{prop} Let $f:M\To E(\sigma)$ be a section with corresponding
function $f^*:G\To E$. Let $X$ be a continuous vector field on $M$
with a horizontal lift $X^*$ to $G$. Then $\nabla_X f$ is the
section which corresponds to the function $X^*f^*=
\mathcal{L}_{X^*}f^*$, where $\mathcal{L}_{X^*}$ denotes the Lie
derivative.
\end{prop}

This same procedure works for all of the bundles in which we have
interest. Note that the two quotient mappings $Q,Q'$ give rise to
different classes of $E$-valued functions, so we get two distinct
connections, $\nabla$ on $E(\sigma)$ and $\nabla'$ on $E(\sigma')$.
It follows from elementary properties of the Lie derivative that
these connections are compatible with our Euclidean structure.

As described in section \ref{subsection:laplacians} these
connections (together with the canonical metric) induce the standard
differential operators: exterior differential, Laplacian, and their
adjoints. We use the standard notations to refer to the objects
associated to $\nabla$ and we use primes to refer to the objects
associated to $\nabla'$. The following result is immediate from the
definitions.
\begin{lem}
  The exterior differential $d$ on
  $\mathcal{E}(\pi)$-valued $p$ forms associated to
  $\nabla$ as in section \ref{subsection:laplacians} is the usual de~Rham
  cohomology coboundary map.
\end{lem}

\subsubsection{An explicit formula.}
\label{subsection:explicitformula}

We now apply Corollary \ref{integrated bw formula} to
$\mathcal{E}(\pi)$ twice, once for each pair of connection and
Laplacian. Taking the difference of the resulting equations and
recalling the notation of section \ref{subsection:laplacians}, we
obtain the following result.

\begin{prop}\label{integrated difference of bw forms} For an
$\mathcal{E}(\pi)$-valued $p$-form $\sigma$ on $M$ we have that
\begin{equation*}
    (\triangle \sigma, \sigma) - (\triangle'\sigma, \sigma) =
    (\nabla^*\nabla\sigma, \sigma) -
    (\left(\nabla'\right)^*\nabla'\sigma,\sigma) +  (
    S(\sigma)-S'(\sigma), \sigma) .
\end{equation*}
\end{prop}

Our goal is now to compute workable expressions for the terms on
the right hand side of this equation in terms of the group $G$ and
the representation $\pi$.

Let $\eta$ be an $\mathcal{E}(\pi)$-valued $p$-form on $M$. We
lift $\eta$ to the corresponding $E$-valued $p$-form $\eta^*$ on
$G$ through the quotient map $Q'$. Thus $\eta^*$ satisfies the
following properties:
\begin{list}
{\bf{P--\Roman{bean}}}{\usecounter{bean}}
    \item $\eta^* \circ L_{\gamma}=(\mathrm{Id}_F\otimes\pi_2(\gamma))\eta^*$
     for $\gamma \in \Gamma$, where
    $L_{\gamma}$ denotes left multiplication,
    \item $i(X)\eta^* = 0$ for $X\in \mathfrak{k}$,
    \item $\eta^*\circ R_k =(\pi_1(k)^{-1}\otimes \mathrm{Id}_{\mathcal{H}})\eta^*$,
    for $k\in K$, where $R_k$ denotes right multiplication, and
    \item $\mathcal{L}_X\eta^*= (-\pi_1(X)\otimes \mathrm{Id}_{\mathcal{H}})\eta^*$
    for $X\in\mathfrak{k}$, where $\mathcal{L}_X$ denotes the Lie
    derivative with respect to $X$.
\end{list}
The fourth property is the linearization of the $K$-equivariance in
property P--III.

As usual, we identify elements of the Lie algebra $\mathfrak{g}$
with left-invariant vector fields on $G$. Also, we fix a basis of
$\mathfrak{g}$ as follows. Let $N=\dim M = \dim \mathfrak{p}$ and
choose a basis $X_1, \ldots, X_N$ of $\mathfrak{p}$ which is
orthonormal with respect to the restriction of the Killing form of
$\mathfrak{g}$ to $\mathfrak{p}$. Let $m = \dim \mathfrak{k}$,
$n=N+m$ and $X_{N+1}, \ldots, X_n$ be an orthonormal basis of
$\mathfrak{k}$ with respect to the negative of the restriction of
the Killing form of $\mathfrak{g}$ to $\mathfrak{k}$.

Let $\omega^1,\ldots, \omega^n$ be the $1$-forms on $G$ which are
dual to the left invariant vector fields $X_1,\ldots, X_n$. These
are determined by the conditions $\omega^j(X_i) = \delta^j_i$. Note
that this means that for $1\leq k,j\leq N$, we have $X_k\omega^j
=0$. We choose our locally symmetric metric on $M$ to be normalized
so that $G{\rightarrow}M$ is a Riemannian submersion for the metric
$\sum_{i=1}^n \omega^i\omega^i$ on $G$. The induced volume form of
$K \subseteq G$ is given by
$\omega^{N+1}\wedge\cdots\wedge\omega^{N+m}$, where the forms are
restricted to the compact Lie group $K$. Also choose a countable
complete orthonormal spanning set $\{v_{\alpha}\}_{\alpha \in A}$
for our Hilbert space $\mathcal{H}$.

We write our $E$-valued $p$-form $\eta^*$ on $G$ as
$$\eta^*= \sum_{\alpha} \eta^{\alpha}\otimes v_{\alpha},$$
where the $v_{\alpha} \in \mathcal{H}$ and the $\eta^{\alpha}$ are
$F$-valued $p$-forms on $G$.

In fact, we identify $\eta^*$ with the system of functions
$\eta^{\alpha}_{i_1,\ldots,i_p}$ on $G$ (with values in $E$)
determined by
$$\eta^*(X_{i_1}, \ldots, X_{i_p}) =
    \sum_{\alpha} \eta^{\alpha}_{i_1,\ldots,i_p}\otimes v_{\alpha}.$$ Note
that we may now write
\begin{equation*}
    \eta^* = \sum_{\alpha\in A} \sum_{i_1<\ldots<i_p}
    \omega^{i_1}\wedge\ldots\wedge\omega^{i_p}\otimes
        \eta^{\alpha}_{i_1,\ldots,i_p}\otimes v_{\alpha}.
\end{equation*}

We compute an expression for the terms on the right hand side of
Corollary \ref{integrated difference of bw forms} in terms of the
functions $\eta^{\alpha}_{i_1,\ldots,i_p}$.

First, we shall need an expression for the inner product of two
forms. Suppose that we have two $\mathcal{E}(\pi)$-valued
$p$-forms $\eta, \theta$ on $M$ with lifts
\begin{align*}
    \eta^*& = \sum_{\alpha\in A} \sum_{i_1<\ldots<i_p}
        \omega^{i_1}\wedge\ldots\wedge\omega^{i_p}\otimes
        \eta^{\alpha}_{i_1,\ldots,i_p}\otimes v_{\alpha}\\
    \theta^*& = \sum_{\beta\in A} \sum_{j_1<\ldots<j_p}
        \omega^{j_1}\wedge\ldots\wedge\omega^{j_p}\otimes
        \theta^{\beta}_{j_1,\ldots,j_p}\otimes v_{\beta}\\
\end{align*}

\noindent Let $c^{-1}$ be the Haar volume of the maximal compact
subgroup $K$ of $G$ and let $D\subset G$ be a fundamental domain
for the action of $\Gamma$ on $G$ by left multiplication. We have
the following modification of a result of
\cite{Matsushima-Murakami}.

\begin{prop}\label{inner product in terms of eta^*}
  For $\mathcal{E}(\pi)$-valued $p$-forms $\eta, \theta$ as above,
  their inner product is
  \begin{equation*}
    (\eta,\theta) = \dfrac{c}{p!}\sum_{\alpha \in A}
    \sum_{i_1,\ldots,i_p = 1}^N \int_D
    (\eta^{\alpha}_{i_1,\ldots,i_p},
    \theta^{\alpha}_{i_1,\ldots,i_p})_F \ dvol_G.
\end{equation*}
\end{prop}

\begin{proof} This argument should be compared with
\cite[Proposition 5.1]{Matsushima-Murakami}.

Recall the quotient mappings $q:G\To M$ and $Q':G\times E \To
E(\sigma)=\mathcal{E}(\pi)$. For $y\in G$, let $Q'_y(v)$ be the
point $Q'(y,v)$ in $\mathcal{E}(\pi)$. Then $v\mapsto Q'_y(v)$ is a
linear isomorphism of $E$ onto the fiber $\mathcal{E}(\pi)_{q(y)}$.
So for each point $y \in G$ and any tangent vectors $Z_1,\ldots,
Z_p$ to $G$ at $y$ we have
$$\eta^*_y(Z_1,\ldots, Z_p) =
{Q'_y}^{-1}\eta_{q(y)}(q_*Z_1,\ldots,q_*Z_p).$$ In particular, this
means that
$$\sum_{\alpha\in A} \eta^{\alpha}_{i_1,\ldots,i_p}\otimes
v_{\alpha} = {Q'_y}^{-1}\eta_{q(y)}(q_*X_{i_1},\ldots,q_*X_{i_p}).$$

Now fix $y \in G$ and let $x=q(y) \in M$. By our choices, the
vectors $$q_*X_1(y), \ldots, q_*X_N(y)$$ form an orthonormal basis
of $T_xM$. In a sufficiently small neighborhood $U$ of $x$, we may
choose an orthonormal basis of $1$-forms $\theta^1,\ldots,\theta^N$
on $M$ such that $$\theta^i_x(q_*X_j(y)) = \delta^i_j.$$

As in section \ref{subsection:laplacians}, we express $\eta$ and
$\theta$ on $U$ as
\begin{eqnarray*}
  \eta & = \sum_{i_1<\ldots<i_p}\omega^{i_1}\wedge\ldots\wedge\omega^{i_p}
   \otimes a_{i_1,\ldots,i_p},\\
  \theta & = \sum_{j_1<\ldots<j_p}\omega^{j_1}\wedge\ldots\wedge\omega^{j_p}
   \otimes b_{j_1,\ldots,j_p}
\end{eqnarray*}
where $a_{i_1,\ldots,i_p}$ and $b_{j_1,\ldots,j_p}$ are sections
of $\mathcal{E}(\pi)$. By the last paragraph, we see that
$$Q'_y\left(\sum_{\alpha\in A} \eta^{\alpha}_{i_1,\ldots,i_p}\otimes
v_{\alpha}\right) = \eta_{q(y)}(q_*X_{i_1},\ldots,q_*X_{i_p}) =
a_{i_1,\ldots,i_p}(x),$$ and
$$Q'_y\left(\sum_{\beta\in A} \theta^{\beta}_{i_1,\ldots,i_p}\otimes
v_{\beta}\right) = \theta_{q(y)}(q_*X_{i_1},\ldots,q_*X_{i_p}) =
b_{j_1,\ldots,j_p}(x).$$

By equation \ref{inner product integrand locally}, we have
$$\eta^t\wedge*\#\theta (x) = \sum_{i_1<\ldots<i_p}
\left(a_{i_1,\ldots,i_p}(x),
b_{j_1,\ldots,j_p}(x)\right)^{\text{can}}(x),$$ and by the
construction of the canonical metric in the fibers of
$\mathcal{E}(\pi)$,
\begin{equation}\label{local integrand lifted}
\begin{split}
  \eta^t\wedge*\#\theta (x) & =
    \sum_{\alpha,\beta}\sum_{i_1<\ldots<i_p}
    \left(\eta^{\alpha}_{i_1,\ldots,i_p}(y),
    \theta^{\beta}_{j_1,\ldots,j_p}(y)\right)_F
    \left(v_{\alpha},v_{\beta}\right)_{\mathcal{H}}\\
 & =\sum_{\alpha}\sum_{i_1<\ldots<i_p}
    \left(\eta^{\alpha}_{i_1,\ldots,i_p}(y),
    \theta^{\alpha}_{j_1,\ldots,j_p}(y)\right)_F .
\end{split}\end{equation}

\noindent The reader should note the disappearance of terms coming
from $\fh$ in the last line of the computation and compare it to
considerations in \cite{Mok}. It is essential for all of our
computations that no terms arise from the ``flat, Hilbertian" part
of either $F{\otimes}\fh$ bundle over $M$.  The quotient map $q:G\To
M$ decomposes as $q=q_1\circ q_2$ where $q_2:G\To \Gamma\backslash
G$ and $q_1:\Gamma\backslash G\To M$. As noted as \cite[Lemma
5.3]{Matsushima-Murakami}, with our conventions, if $f$ is a
continuous function on $M$, then $\int_M f\ dvol_M = c
\int_{\Gamma\backslash G} f\circ q_1 dvol_{\Gamma\backslash G}$. By
lifting one more step up through $q_2$ we obtain that $\int_M f\
dvol_M = c \int_D f\circ q \ dvol_{G}.$ A simple application of this
last equation to equation \ref{local integrand lifted} completes the
proof of the Proposition.\end{proof}

We now compute the terms on the right hand side of \ref{integrated
difference of bw forms} in terms of the functions
$\eta^{\alpha}_{i_1,\ldots,i_p}$. We begin by working with the
connection in $E(\sigma')$.

Fix $y$ in $G$, and let $x= q(y)=\Gamma y K \in M$. The vectors
$X_1(y),\ldots,X_N(y)$ project to an orthonormal basis
$\{q_*X_i(y)\}$ of $T_xM$. On a sufficiently small neighborhood $U$
of $x$ in $M$, we extend this basis to a synchronous framing
$\{Y_i\}$ over $U$ such that $Y_i(x) = q_*X_i(y)$,
$\nabla_{Y_j}Y_i(x)=0$, and $(Y_j,Y_k) = \delta_{jk}$ on $U$. Recall
that such a framing is constructed by parallel transport of the
frame $\{q_*X_i(y)\}$ radially outward from $x$ along geodesics.

Let $\theta^i$ be the basis for ordinary $1$-forms on $U$ which is
dual to the framing $\{Y_i\}$. That is, at each point of $U$,
$\theta^i(Y_j) = \delta^i_j$. Note that this means that the
$\theta^i$'s form an orthonormal basis for the $1$-forms on $U$. As
above, we decompose $\eta$ in terms of this basis.
$$\eta = \sum_{i_1<\ldots<i_p} a_{i_1,\ldots,i_p}\otimes
\theta^{i_1}\wedge\ldots\wedge\theta^{i_p}.$$

So we compute
\begin{multline*}
    \nabla_{Y_k} \eta = \sum_{i_1<\ldots<i_p}
    \nabla_{Y_k}a_{i_1,\ldots,i_p} \otimes
    \theta^{i_1}\wedge\ldots\wedge\theta^{i_p}\\
    + \sum_{i_1<\ldots<i_p}\sum_{s=1}^p a_{i_1,\ldots,i_p} \otimes
    \theta^{i_1}\wedge\ldots\wedge\nabla_{Y_k}\theta^{i_s}
    \wedge\ldots\wedge\theta^{i_p}.
\end{multline*}
Note that because we are working in a synchronous frame, we have
\begin{equation*}
  \nabla_{Y_k}\theta^i (x) = 0 ,\quad  \text{for all $i$ and $k$}.
\end{equation*}
Hence the second sum vanishes and we need to compute $\nabla_{Y_k}
a_{i_1,\ldots,i_p}(x)$.

Since $a_{i_1,\ldots,i_p}$ is a section of $E(\sigma)$ we may use
the Lie derivative description from section \ref{connections}.
Choose horizontal lifts $Y^*_k$ of the vector fields $Y_k$. By
definition of the lifting procedure, the lift $a^*_{i_1,\ldots,i_p}$
of $a_{i_1,\ldots,i_p}$ is $$a_{i_1,\ldots,i_p}^* =
{Q'_y}^{-1}(\eta(Y_{i_1},\ldots,Y_{i_p})) =
\eta^*(Y^*_{i_1},\ldots,Y^*_{i_p}).$$ So that
\begin{equation*}
    a_{i_1,\ldots,i_p}^* =\sum_{\alpha \in A}
    \sum_{j_1<\ldots<j_p}\eta^{\alpha}_{j_1,\ldots,j_p}\otimes
    v_{\alpha}\cdot
    \omega^{j_1}(Y^*_{i_1})\cdots\omega^{j_p}(Y^*_{i_p})
\end{equation*}
Hence the function $G\To E$ corresponding to
$\nabla_{Y_k}a_{i_1,\ldots,i_p}$ is
\begin{equation*}
\begin{split}
    L(y) &= Y^*_k(y)a^*_{i_1,\ldots,i_p}(y)\\ &= \sum_{\alpha\in
    A}\sum_{j_1<\ldots<j_p} Y^*_k\eta^{\alpha}_{j_1,\ldots,j_p}\otimes
    v_{\alpha}\cdot
    \omega^{j_1}(Y^*_{i_1})\cdots\omega^{j_p}(Y^*_{i_p})\\
    &+\sum_{\alpha\in A}\sum_{j_1<\ldots<j_p} \sum_{s=1}^p
    \eta^{\alpha}_{j_1,\ldots,j_p}\otimes v_{\alpha}\cdot
    \omega^{j_1}(Y^*_{i_1})\cdots\{Y^*_k(\omega^{j_s}(Y^*_{i_s}))\}
    \cdots\omega^{j_p}(Y^*_{i_p})
\end{split}
\end{equation*}
Now, by the definition of the Lie derivative, to compute
$Y^*_k(\omega^{j_s}(Y^*_{i_s}))$ we only need to know
$\omega^{j_s}(Y^*_{i_s})$ along a curve $\gamma(t)$ in $G$ which
passes through $\gamma(0)=y$ in the direction $\gamma'(0)=X_k(y)$.
The assignment $t\mapsto y\exp(tX_k(e))$ is easily checked to be
such a curve. In fact, this curve is horizontal since $\gamma'(t)
\in H_{\gamma(t)}$ for all $t$. So by the definition of parallel
translation, we see that
$$Y_i^*(\gamma(t)) = X_i(\gamma(t)).$$ This gives that, at a point
$\gamma(t)$, $\omega^{j_r}(Y^*_{i_r}) = \delta^{j_r}_{i_r}$ and
hence along $\gamma(t)$ we have
$$Y^*_k(\omega^{j_s}(Y^*_{i_s})) = 0.$$
Using this and that $Y^*_i(y)=X_i(y)$, we see that
\begin{equation*}
  L(y) = X_k(y)a^*_{i_1,\ldots,i_p}(y) = \sum_{\alpha\in
    A} X_k\eta^{\alpha}_{i_1,\ldots,i_p}\otimes
    v_{\alpha}.
\end{equation*}
From this argument, we deduce the following lemma.

\begin{lem}\label{nabla'}
  $\left(\nabla'_{Y_k}\eta \right)^{\alpha}_{i_1,\ldots,i_p}(y) =
  X_k\eta^{\alpha}_{i_1,\ldots,i_p}(y)$
\end{lem}

Our next step is to compute the adjoint of $\nabla'_{Y_k}$ at $y$.
We need the following result of Weil, which is a form of integration
by parts.
\begin{prop}[Weil, \cite{Weil-discrete-subgroups-II}]
    Let $f$ be a smooth function on $\Gamma\backslash G$, then
    $$\int_{\Gamma\backslash G} X_i f \ dvol_{\Gamma\backslash G}
    =0.$$
\end{prop}
\noindent Lifting this result to the fundamental domain $D$ for the
$\Gamma$-action on $G$ we have:
\begin{equation}
\label{equation:weillifted} \int_D X_i f \ dvol_G
    =0
\end{equation}
\noindent for every  smooth $\Gamma$-invariant function $f$ on $G$.

For any two $p$-forms $\eta, \theta$ we have
\begin{equation*}
((\nabla')_{Y_k}^*\eta,\theta) =\dfrac{c}{p!}\sum_{\alpha \in A}
    \sum_{i_1,\ldots,i_p = 1}^N \int_D
    \left(((\nabla')_{Y_k}^*\eta)^{\alpha}_{i_1,\ldots,i_p},
    \theta^{\alpha}_{i_1,\ldots,i_p}\right)_F \ dvol_G,
\end{equation*}
\noindent and by Lemma \ref{nabla'}, equation
\ref{equation:weillifted}, and the compatibility of the connection
and the fiberwise metric, we see that
\begin{equation*}
\begin{split}
 (\eta, \nabla'_{Y_k}\theta)&=\dfrac{c}{p!}\sum_{\alpha \in A}
    \sum_{i_1,\ldots,i_p = 1}^N \int_D
    \left(\eta^{\alpha}_{i_1,\ldots,i_p},
    X_k\theta^{\alpha}_{i_1,\ldots,i_p}\right)_F \ dvol_G  \\
    &= \dfrac{c}{p!}\sum_{\alpha \in A}
    \sum_{i_1,\ldots,i_p = 1}^N \int_D
    \left(-X_k\eta^{\alpha}_{i_1,\ldots,i_p},
    \theta^{\alpha}_{i_1,\ldots,i_p}\right)_F \ dvol_G.
\end{split}
\end{equation*}

\noindent By the equality $((\nabla')_{Y_k}^*\eta,\theta)=(\eta,
\nabla'_{Y_k}\theta)$ defining the adjoint, the following result
holds.

\begin{lem}\label{nabla' adjoint}
 $\left((\nabla')^*_{Y_k}\eta \right)^{\alpha}_{i_1,\ldots,i_p}(y) =
  -X_k\eta^{\alpha}_{i_1,\ldots,i_p}(y)$.
\end{lem}

Using these results together, we can now find an expression for one
of the terms of the right hand side of Proposition \ref{integrated
difference of bw forms}.

\begin{prop}\label{K covariant laplacian}
For an $\mathcal{E}(\pi)$valued $p$-form $\eta$,
  $$((\nabla')^*\nabla'\eta)^{\alpha}_{i_1,\ldots,i_p} = -\sum_{k=1}^N
    X_k^2\eta^{\alpha}_{i_1,\ldots,i_p},$$
and hence
  $$((\nabla')^*\nabla'\eta, \eta) = \dfrac{c}{p!} \sum_{\alpha\in A}
  \sum_{i_1,\ldots,i_p=1}^N \int_{D}\left(-\sum_{k=1}^N
    X_k^2 \eta^{\alpha}_{i_1,\ldots,i_p},\
    \eta^{\alpha}_{i_1,\ldots,i_p}\right)_F \ dvol_{G}.$$
\end{prop}

\begin{proof} The second equality follows directly from the first
and Proposition \ref{inner product in terms of eta^*}. For the
first note that at $x$ we have
\begin{equation*}
    {\nabla'}^*\nabla' \eta (x) = -\sum_k {\nabla'}^*_{Y_k}\nabla'_{Y_k}\eta
    (x).
\end{equation*}
\noindent This is a standard computation---see for example
\cite[page 73]{Rosenberg} for the case of ordinary forms, the
proof for vector bundle valued forms is analogous. The proposition
now follows by an application of Lemmas \ref{nabla'} and
\ref{nabla' adjoint}.\end{proof}

We make similar computations for the connection $\nabla$ in
$E(\sigma)$. As much is the same, we highlight differences and
state results.

The main difference is that to use the Lie derivative description
of $\nabla$ one must instead lift into the bundle $E(\sigma)$,
while the functions $\eta^{\alpha}_{i_1,\ldots,i_p}$ are lifts
into $E(\sigma')$. We remedy this by using the isomorphism between
the two bundles given by the map $A:G\times E \To G\times E$
defined in Lemma \ref{lifted identification of bundles}.
Horizontal lifts of vectors tangent to $M$ are performed in
exactly the same way, so no change is required on this point.

We pull functions back by $A$, take their appropriate Lie
derivatives and push forward by $A$. Inspecting the definition of
$A$, we see that the relevant functions in the $E(\sigma)$ lift are
then $y\mapsto \pi_1(y)\eta^{\alpha}_{i_1,\ldots,i_p}(y)$. We then
take the Lie derivative
\begin{equation*}
\begin{split}
    X_k (\pi_1(y)\eta^{\alpha}_{i_1,\ldots,i_p}(y))&
    = (\mathcal{L}_{X_k}(\pi_1))(y)\eta^{\alpha}_{i_1,\ldots,i_p}(y)
    + \pi_1(y) X_k\eta^{\alpha}_{i_1,\ldots,i_p}(y) \\
    &= \pi_1(y)\pi_1(X_k)\eta^{\alpha}_{i_1,\ldots,i_p}(y)
    + \pi_1(y) X_k\eta^{\alpha}_{i_1,\ldots,i_p}(y),
\end{split}
\end{equation*}
where $(\mathcal{L}_{X_k}(\pi_1))$ is the Lie derivative of $\pi$
along $X_k$ and $\pi_1(X_k)$ denotes the image of $X_k\in
\mathfrak{p}$ by the differential of the homomorphism $\pi_1$ at
the identity. Applying $A$ again, we deduce the following formula.
\begin{lem}\label{nabla}
  $\left(\nabla_{Y_k}\eta \right)^{\alpha}_{i_1,\ldots,i_p}(y) =
  (X_k+\pi_1(X_k))\eta^{\alpha}_{i_1,\ldots,i_p}(y)$
\end{lem}

The adjoint of $\nabla_{Y_k}$ is computed as before, using the
additional observation that for $X_k \in \mathfrak{p}$, the linear
map $\pi_1(X_k)$ is symmetric with respect to $(\ ,)_F$. The
corresponding formula is:
\begin{lem}\label{nabla adjoint}
 $\left(\nabla^*_{Y_k}\eta \right)^{\alpha}_{i_1,\ldots,i_p}(y) =
  -(X_k-\pi_1(X_k))\eta^{\alpha}_{i_1,\ldots,i_p}(y)$
\end{lem}

Combining the last two lemmas, we have:

\begin{lem}\label{Gamma covariant laplacian}
For an $\mathcal{E}(\pi)$-valued $p$-form $\eta$,
  $$\left(\nabla^*\nabla\eta\right)^{\alpha}_{i_1,\ldots,i_p} = -\sum_{k=1}^N
    (X_k^2 -\pi(X_k)^2)\eta^{\alpha}_{i_1,\ldots,i_p},$$
and hence
\begin{equation*}
  \left(\nabla^*\nabla\eta, \eta\right) =
  \dfrac{c}{p!} \sum_{\alpha\in A}
  \sum_{i_1,\ldots,i_p=1}^N \int_{D}\left(-\sum_{k=1}^N
    (X_k^2 -\pi(X_k)^2)\eta^{\alpha}_{i_1,\ldots,i_p},\
    \eta^{\alpha}_{i_1,\ldots,i_p}\right)_F \ dvol_{G}.
    \end{equation*}
\end{lem}

\begin{proof} This lemma is proved in the same way as Proposition \ref{K
covariant laplacian}, using Lemmas \ref{nabla} and \ref{nabla
adjoint} in place of Lemmas \ref{nabla'} and \ref{nabla'
adjoint}.\end{proof}

Combining Lemmas \ref{K covariant laplacian} and \ref{Gamma
covariant laplacian} yields:

\begin{lem}
\begin{multline*}  (\nabla^*\nabla\eta, \eta) -
    (\left(\nabla'\right)^*\nabla'\eta,\eta) = \\\dfrac{c}{p!} \sum_{\alpha\in A}
  \sum_{i_1,\ldots,i_p=1}^N \int_{D}\left(\sum_{k=1}^N
    \pi(X_k)^2 \eta^{\alpha}_{i_1,\ldots,i_p},\
    \eta^{\alpha}_{i_1,\ldots,i_p}\right)_F \
    dvol_{G}.\end{multline*}
\end{lem}

It remains only to compute the curvature term from Proposition
\ref{integrated difference of bw forms}. By definition, the
curvature terms are zero if $p=0$, so we assume $p\geq1$. In the
definition of the Ricci tensor $S$ of a vector bundle, one must
use the curvature tensor $R^V$ of the connected bundle. The
Leibniz property for $R^V$ for vector fields $X,Y,Y_1,\dots, Y_p$
is
\begin{equation}
\begin{split}
(R^V(X,Y)\eta)(Y_1,\ldots,Y_p) &= R^V(X,Y)(\eta(Y_1,\ldots,Y_p))
\\ &\qquad- \sum_{i=1}^p \eta(Y_1,\ldots, R^N(X,Y)Y_i, \ldots, Y_p).
\end{split}
\end{equation}
\noindent Note that the sum in the second term is independent of the
connection in the bundle: it depends only on the structure of the
Riemannian manifold. Hence, the difference of the corresponding
terms in Proposition \ref{integrated difference of bw forms} is
zero. Also, the $\Gamma$-bundle connection is flat by construction,
so the first term in this case is zero. This leaves us only one
term. For a point $x\in M$ and vectors $Y_1,\dots, Y_p$ at $x$ and
an orthnormal basis $e_i$ of $T_xM$,
\begin{equation}
\begin{split}
(S(\eta)_x-S'(\eta)_x &)(Y_1, \ldots, Y_p)=\\
    &\sum_{k=1}^p\sum_{s=1}^n (-1)^{k+1}
        R'(e_s,Y_k)(\eta(e_s,Y_1, \ldots, \widehat{Y_k}, \ldots,
        Y_p)).
\end{split}
\end{equation}

\noindent Recall that the $Y_i$'s form an orthonormal basis in a
neighborhood $U$ of $x$. By the definition of curvature as an
iterated covariant derivative and Lemma \ref{nabla'}, we see that
\begin{equation}
\begin{split}
    (S(\eta)-S'(\eta))^\alpha_{i_1,\ldots,i_p}& =
    \sum_{k=1}^p\sum_{s=1}^N(-1)^{k+1} (-X_sX_{i_k}+
    X_{i_k}X_s)\eta^{\alpha}_{s,i_1,\ldots,\hat{i_k},\ldots,i_p}\\
    &=\sum_{k=1}^p \sum_{s=1}^N (-1)^{k+1} [X_{i_k}, X_s]
    \eta^{\alpha}_{s,i_1,\ldots,\hat{i_k},\ldots,i_p}
\end{split}
\end{equation}
By property P-IV of our forms, since $[X_h,X_{i_k}]\in
\mathfrak{k}$, we obtain
\begin{equation}
    (S(\eta)-S'(\eta))^\alpha_{i_1,\ldots,i_p}= \sum_{k=1}^p \sum_{s=1}^N (-1)^{k}
    \pi_1([X_{i_k}, X_s])\eta^{\alpha}_{s,i_1,\ldots,\widehat{i_k},\ldots,i_p}
\end{equation}

\noindent Applying the formula for inner products of forms, we
obtain
\begin{equation}
\left(S(\eta)-S'(\eta), \eta\right) = \dfrac{c}{p!} \sum_{\alpha \in
A} \sum_{i_1,\ldots,i_p=1}^N
    \int_{D} \left( (S(\eta)-S'(\eta))^\alpha_{i_1,\ldots,i_p},
    \eta^{\alpha}_{i_1,\dots,i_{p}}\right)_F\ dvol_G.
\end{equation}

Combining the last two equations and doing a small bit of
re-indexing yields the last piece of Proposition \ref{integrated
difference of bw forms}.

\begin{lem}\label{curvature term} The final term of the right hand
side of Proposition \ref{integrated difference of bw forms} is
\begin{multline}
\left( S(\eta)-S'(\eta), \eta \right) \\
=\dfrac{c}{(p-1)!}\sum_{\alpha\in A}\sum_{h,i_1,\ldots,i_{p-1}=1}^N
    \int_{D} \left( \sum_{s=1}^N
    \pi_1([X_h,X_s])\eta^{\alpha}_{s,i_1,\dots,i_{p-1}},
    \eta^{\alpha}_{h,i_1,\dots,i_{p-1}} \right)_F\ dvol_G.
\end{multline}
\end{lem}

Assembling the results of Corollary \ref{integrated bw formula}
and Lemmas \ref{K covariant laplacian}, \ref{Gamma covariant
laplacian} and \ref{curvature term}, we deduce our variant of the
Bochner formula.
\begin{prop}\label{Our Bochner formula} Let $p\geq 1$.
Let $\eta$ be an $\mathcal{E}(\pi)$-valued $p$-form on $M$ with
$E(\sigma')$ lift $$\eta^* = \sum_{\alpha \in A}
\sum_{i_1<\ldots<i_p} \eta^{\alpha}_{i_1,\ldots,i_p}\otimes
v_{\alpha}\otimes \omega^{i_1}\wedge\ldots\wedge \omega^{i_p}$$ as
above. Then
\begin{multline}
( \triangle \eta, \eta) -( \triangle'\eta, \eta)  =\\
    \dfrac{c}{p!}\sum_{\alpha \in A}\sum_{i_1,\ldots,i_p=1}^N \int_{D}
    \left(\sum_{k=1}^N\pi_1(X_k)^2\eta^{\alpha}_{i_1,\ldots,i_p},
    \eta^{\alpha}_{i_1,\ldots,i_p}\right)_F\ dvol_G  \\
 + \dfrac{c}{(p-1)!}\sum_{\alpha\in A}\sum_{h,i_1,\ldots,i_{p-1}=1}^N
    \int_{D} \left( \sum_{k=1}^N
    \pi_1([X_h,X_k])\eta^{\alpha}_{k,i_1,\dots,i_{p-1}},
    \eta^{\alpha}_{h,i_1,\dots,i_{p-1}}\right)_F\ dvol_G.
\end{multline}
\end{prop}

\subsubsection{Estimate on the Laplacian}\label{subsection:endgame}

In this section, we use results of Raghunathan to complete the proof
of Theorems \ref{thm:1-vanishing} and \ref{estimate on forms}. In
fact, we obtain more general statements in both cases, which depend
upon more detailed knowledge of the representation.

First, it will be useful to rearrange the result of Proposition
\ref{Our Bochner formula} slightly. We have the following equivalent
formulation.
\begin{equation}\label{compact bochner formula}
(\triangle\eta, \eta) - (\triangle'\eta, \eta)  = \\
\dfrac{c}{(p-1)!}\sum_{\alpha\in A}\int_{D}\sum_{i_1,\ldots,i_p=1}^N
\left(T^p_{\pi_1}\eta^{\alpha}_{i_1,\ldots,i_p},
\eta^{\alpha}_{i_1,\ldots,i_p}\right)_F\ dvol_G,
\end{equation}
where the function $T^p_{\pi_1}\eta^{\alpha}_{i_1,\ldots,i_p}:G\To
F$ is defined by
\begin{equation}\label{endomorphism}
  T^p_{\pi_1}\eta^{\alpha}_{i_1,\ldots,i_p} =\frac{1}{p}\sum_{k=1}^N\pi_1(X_k)^2\eta^{\alpha}_{i_1,\ldots,i_p}
+\sum_{k=1}^N
\pi_1([X_{i_1},X_k])\eta^{\alpha}_{k,i_1,\dots,i_{p-1}}.
\end{equation}

We introduce some terminology which will bring us closer to the
language in
\cite{Raghunathan-1-vanishing,Raghunathan-highervanishing,Raghunathan-correction}.

Given fixed $p\geq 1$, $y\in G$, $\alpha\in A$ we shall consider
$\theta = \eta^{\alpha}_y$ as an element of the finite dimensional
vector space $R_p := \bigwedge^p\mathfrak{p}^*\otimes F$. We define
an inner product $I_p$ on $R_p$ as follows: for $\theta^1,\theta^2
\in R_p$, put
$$I_p(\theta^1,\theta^2)= \sum_{i_1,\ldots,i_p=1}^N
\left(\theta^1(X_{i_1},\ldots,X_{i_p}),
\theta^2(X_{i_1},\ldots,X_{i_p})\right)_F.$$ Also, the mapping
$T^p_{\pi_1}:R_p\To R_p$ defined by
\begin{equation*}
T^p_{\pi_1}\theta (Y_1,\ldots,Y_p)=
\sum_{k=1}^N\left\{\frac{1}{p}\pi_1(X_k)^2\theta (Y_1,\ldots,Y_p)
+\pi_1([Y_1,X_k])\theta (X_k,Y_2,\ldots,Y_p)\right\}
\end{equation*}
is a linear endomorphism of $R_p$. It is not difficult to check that
this is symmetric with respect to $I_p$, so all of its eigenvalues
are real.

\begin{prop} To check the hypothesis of Proposition
\ref{sufficient condition for vanishing}, it suffices to show that
the quadratic form
\begin{equation*}
\theta \To I_p(T^p_{\pi_1}\theta, \theta)
\end{equation*}
on $R_p$ is positive definite.
\end{prop}

\begin{proof}
Since $R_p$ is finite dimensional, positive-definiteness is
equivalent to having all the eigenvalues of $T^p_{\pi_1}$
positive. Then one can choose $C$ to be the minimal eigenvalue of
$T^p_{\pi_1}$ and obtain $I_p(T^p_{\pi_1}\theta, \theta) \geq C
\cdot I_p(\theta,\theta)$ for every $\theta \in R_p$. It is
important to note that $C$ depends only upon $p$ and $\pi_1$.

For each $y\in G$ and for each $\alpha \in A$, we think of
$\eta^{\alpha}_y$ as an element of $R_p$. Using the definitions of
$T^p_{\pi_1}$ and $I_p$ along with equations \ref{compact bochner
formula} and \ref{endomorphism} and Proposition \ref{inner product
in terms of eta^*}, we deduce that
\begin{equation*}
\begin{split}
(\triangle\eta, \eta) &\geq
(\triangle\eta, \eta) -  (\triangle'\eta, \eta)  \\
&=\dfrac{c}{(p-1)!}\sum_{\alpha\in
A}\int_{D}\sum_{i_1,\ldots,i_p=1}^N
\left(T^p_{\pi_1}\eta^{\alpha}_{i_1,\ldots,i_p}(y),
\eta^{\alpha}_{i_1,\ldots,i_p}(y)\right)_F\ dvol_G\\
&=\dfrac{c}{(p-1)!}\sum_{\alpha\in A}\int_{D}
I_p(T^p_{\pi_1}\eta^{\alpha}_y, \eta^{\alpha}_y)\ dvol_G \\
& \geq C\cdot\dfrac{c}{(p-1)!}\sum_{\alpha\in A}\int_{D}
I_p(\eta^{\alpha}_y, \eta^{\alpha}_y)\ dvol_G \\
&=C\cdot p \cdot(\eta,\eta).
\end{split}
\end{equation*}
Relabel $C p$ as $C$ to complete the proof. \end{proof}

Let $\mathfrak{g}_{\mathbb{C}}$ be the complexification of
$\mathfrak{g}$. For a representation $\pi$ of $\mathfrak{g}$ on a
real vector space, $F$, we let $\pi_{\mathbb{C}}$ be the extension
of $\pi$ to a representation of $\mathfrak{g}_{\mathbb{C}}$ on
$F_{\mathbb{C}} = F\otimes_{\mathbb{R}}\mathbb{C}$. Denote by
$\Lambda_{\pi}$ the highest weight of $\pi_{\mathbb{C}}$.

If $\mathfrak{g}=\mathfrak{so}_0(n,1)$ or $\mathfrak{su}(n,1)$ we
denote by $\tau^N$ the natural representation on $\mathbb{R}^{n+1}$
or $\mathbb{R}^{2(n+1)}$, respectively, and let $\mu_N$ be the
highest weight of $\tau^N_{\mathbb{C}}$.

\begin{thm}[{Raghunathan \cite[Theorem 1', page 106]{Raghunathan-1-vanishing}}]
Let $\mathfrak{g}$ be a
semisimple Lie algebra and $\pi$ a non-trivial irreducible
representation on a real vector space $F$. Let $\mathfrak{g}= \sum_i
\mathfrak{g}^i$ be the decomposition of $\mathfrak{g}$ into simple
factors and let $\pi^i = \pi_{\mid_{\mathfrak{g}^i}}$. Suppose that
the following conditions hold:
\begin{itemize}
\item There is an $i$ such that $\mathfrak{g}^i$ is non-compact and
$\pi^i$ is non-trivial.
\item No pair $(\mathfrak{g}^i, \Lambda_{\pi^i})$ is of the form
$(\mathfrak{so}_0(n,1), m\cdot \mu_N)$ or $(\mathfrak{su}(n,1),
m\cdot \mu_N)$, where $m$ is an integer.
\end{itemize}
Then $T^1_{\pi}$ is positive definite.
\end{thm}

Combining this with the results above yields the following
strengthening of Theorem \ref{thm:1-vanishing}.

\begin{thm} \label{full H^1 theorem real case}
Let $G$ be a connected, semisimple real Lie group with Lie algebra
$\mathfrak{g}$. Let $\mathfrak{g}=\sum_i \mathfrak{g}^i$ be the
decomposition of $\mathfrak{g}$ into simple factors, and let
$\mathfrak{g}^i_{\mathbb{C}}$ be the complexification of
$\mathfrak{g}^i$.

Let $\pi_1:G\To \mathrm{GL}(F)$ be a finite dimensional
representation on a real vector space, $F$. Denote by $\pi_1$ also
the induced representation of $\mathfrak{g}$ on $F$. Let $\pi^i =
{\pi_1}_{\mid_{\mathfrak{g}^i}}$, and let $\Lambda_i$ be the highest
weight of the extension of $\pi^i$ to a representation of
$\mathfrak{g}^i_{\mathbb{C}}$ on $F_{\mathbb{C}} =
F\otimes_{\mathbb{R}} \mathbb{C}$.

Let $\Gamma$ be an irreducible cocompact lattice in $G$, and let
$\pi_2$ be a non-trivial unitary representation of $\Gamma$ on a
separable Hilbert space $\mathcal{H}$. Form the representation
$\pi = {\pi_1}_{\mid_{\Gamma}} \otimes \pi_2$.

Suppose that
\begin{itemize}
    \item There is an $i$ such that $\mathfrak{g}^i$ is non-compact and
        $\pi^i$ is non-trivial.
    \item No pair $(\mathfrak{g}^i, \Lambda_{\pi^i})$ is of the form
        $(\mathfrak{so}_0(n,1), m\cdot \mu_N)$ or $(\mathfrak{su}(n,1),
        m\cdot \mu_N)$ for an integer $m$, where $\mu_N$ denotes the
        highest weight of the complexification of the natural representation
        on $\mathbb{R}^{n+1}$ or $\mathbb{R}^{2(n+1)}$, respectively.
\end{itemize}

Then $H^1(\Gamma,\pi) = 0$.
\end{thm}

Note that in the course of the proof, we have also proved the
following extension of Theorem \ref{estimate on forms}.

\begin{thm}\label{full estimate on 1 forms}
Assume the hypotheses and notation of Theorem \ref{full H^1
theorem real case}. Also assume that $\Gamma$ is torsion free.
Form the locally symmetric Riemannian manifold $M =
\Gamma\backslash G/K$, where $K$ is a maximal compact subgroup of
$G$, and form the vector bundle $\mathcal{E}(\pi)$ over $M$
associated to the principal $\Gamma$-bundle $G/K\To M$ and the
representation $\pi$ of $\Gamma$ on $F\otimes \mathcal{H}$.

Then there exists a positive constant $C$ such that for any smooth
$\mathcal{E}(\rho)$-valued differential $1$-form $\eta$ on $M$ we
have
\begin{equation*}
    (\triangle \eta, \eta) \geq  C(\eta, \eta).
\end{equation*}
\end{thm}

\noindent{\bf On Pointwise Formulas:} Though we do not need it in
our applications, we point out that the computations above yield
pointwise formulas as well.  The main point is to subtract the two
versions of equation \ref{equation:pointwisebochnerformula} rather
than the two versions of equation \ref{equation:integratedbw}. The
first point is that the left hand sides of equation
\ref{equation:pointwisebochnerformula} cancel, since the
Laplacians on functions defined by our two connections are equal.
Once one observes that, it is easy to see that all of the
computations in \S\ref{subsection:laplacians} and
\S\ref{subsection:endgame} yield pointwise formulas.  This is
explicit in the statements of most of the lemmas above and easy to
check in the remaining cases.

\subsubsection{Vanishing in higher degrees}
\label{subsection:higherdegrees}

 We obtain results for vanishing in higher degrees in the
case where $\pi_1$ is a representation on a complex vector space $F$
using the computations in
\cite{Raghunathan-highervanishing,Raghunathan-correction}. Note that
none of the development in sections \S\ref{subsection:cohomology
change} through \S\ref{subsection:endgame} depends upon whether $F$
is a real or complex space. The only change is that $R_p$ is a
complex vector space and $I_p$ is a Hermitian product.

To state the result we introduce some more notation. Let $\theta$
denote the Cartan involution of $\mathfrak{g}$ associated to
$\mathfrak{k}$, as well as its extension to
$\mathfrak{g}_{\mathbb{C}}$. Let $\mathfrak{p}_{\mathbb{C}}$ be
the orthogonal complement of $\mathfrak{k}_{\mathbb{C}}$ with
respect to the Killing form $\kappa$ of
$\mathfrak{g}_{\mathbb{C}}$. Let $\mathfrak{h_k}$ be a Cartan
subalgebra of $\mathfrak{k}$, and
$\mathfrak{h}\supset\mathfrak{h_k}$ a Cartan subalgebra of
$\mathfrak{g}$. Let $\mathfrak{h}_{\mathbb{C}}$ (respectively
$\mathfrak{h_k}_{\mathbb{C}}$) be the $\mathbb{C}$-span of
$\mathfrak{h}$ (respectively $\mathfrak{h_k}$).

Let $\Delta$ be the root system of $\mathfrak{g}_{\mathbb{C}}$ with
respect to $\mathfrak{h}_{\mathbb{C}}$. For $\alpha\in \Delta$, let
$H_{\alpha}$ be the unique element of $\mathfrak{h}_{\mathbb{C}}$
such that $\kappa(H, H_{\alpha})= \alpha(H)$ for all $H\in
\mathfrak{H}_{\mathbb{C}}$. Let $\mathfrak{h}^* = \sum_{\alpha\in
\Delta} \mathbb{R}H_{\alpha} = i \mathfrak{h_k}_{\mathbb{C}}\oplus
(\mathfrak{h}_{\mathbb{C}}\cup\mathfrak{p}_{\mathbb{C}})$. We
partition the set of roots $\Delta$ into three subsets. Let
\begin{eqnarray*}
    A&=\{\alpha\in\Delta \mid \theta(E_{\alpha}) = E_{\alpha}\} \\
    B&=\{\alpha\in\Delta \mid \theta(\alpha) \neq \alpha\} \\
    C&=\{\alpha\in\Delta \mid \theta(E_{\alpha}) = -E_{\alpha}\},
\end{eqnarray*}
where $E_{\alpha}$ denotes a root vector corresponding to $\alpha$.
Then we have that $\Delta = A\cup B\cup C$, $\theta$ stabilizes the
sets $\mathfrak{h}_{\mathbb{C}}$, $A\cup C$ and $B$, and
$\theta(\alpha)=\alpha$ if $\alpha \in A\cup C$.

Following Raghunathan, we say that an ordering on the real dual of
$\mathfrak{h}^*$ is \emph{admissible} when it is described by the
following process: Let $H_1,\ldots, H_m$ be an orthonormal basis of
$\mathfrak{h}^*$ chosen so that the initial elements $H_1, \ldots,
H_l$ form a basis of $i\mathfrak{h_k}$, and declare that $\alpha$ in
the real dual of $\mathfrak{h}^*$ is positive if the first
nonvanishing $\alpha{H_i}$ is positive. If $O$ is an admissible
order, then for a finite subset $E$ of the dual of $\mathfrak{h}^*$
we shall write $E^+(O)$ for the subset of positive elements in $E$.

Now, for an irreducible representation $\pi$ of $G$, let
$\Lambda_{\pi}(O)$ denote the highest weight of $\pi$ with respect
to the order $O$. Let $\Sigma_2 (O) = C^+(O) \cup \{\alpha \in
B^+(O) \mid \alpha > \theta(\alpha) \}$ and $\Sigma_{\pi}(O)=
\{\alpha \in \Sigma_2(O) \mid \kappa(\Lambda_{\pi}(O), \alpha)\neq
0\}$. With all of this notation, we can now state the following
result.
\begin{thm} [Raghunathan \cite{Raghunathan-correction,
Raghunathan-1-vanishing}] Let $\pi$ be an irreducible finite
dimensional representation of $G$. Then if $\Sigma_{\pi}(O)$
contains (strictly) more than q elements for every admissible $O$,
then the Hermitian quadratic form  $\eta \mapsto I_p(T^p_{\pi}\eta,
\eta)$ is positive definite for $p\leq q$.\end{thm}

The same line of argument as above allows us to deduce the following
result.
\begin{thm}
  Assume the notations and hypotheses of Theorems
  \ref{full H^1 theorem real case} and \ref{full estimate on 1 forms},
  with the exceptions that $F$ is a complex vector space and $\pi_1$
  is a complex linear representation, and $\Gamma$ is torsion free.
  If Raghunathan's condition on admissible orders is satisfied, then for each
  $1\leq p\leq q$ there is a constant
  $C_p$ such that every $\mathcal{E}(\rho)$-valued $p$-form $\eta$
  on $M$ satisfies
  \begin{equation*}
        (\triangle \eta, \eta) \geq  C_p(\eta, \eta).
  \end{equation*}
Furthermore, $H^p(\Gamma,\rho)=0$ for $1\leq p\leq q$. This last
assertion holds even if $\Gamma$ contains elements of finite order.
\end{thm}

\subsection{Direct Computation of the flat
Laplacian}\label{subsection:flatlaplacian}

As it will be useful for the discussion to follow in \S
\ref{section:gkcohomology}, we now give a more direct computation of
the Laplacian $\triangle$ associated to the flat connection
$\nabla$. The main point is that we may realize the action of the
Laplacian as the action of the Casimir operator of our Lie algebra
on the family of functions $\eta^{\alpha}_{i_1,\ldots,i_p}$. To do
so, we continue with the notations and setup of this section. In
particular, we make use of Lemma \ref{nabla}.

Recall that the exterior differential $d$ may be defined as the
anti-symmetrization of the connection. So for our
$\mathcal{E}(\rho)$-valued $p$ form $\eta$ on $M$ and vector fields
$Z_1, \ldots, Z_{p+1}$ on $M$ we have
\begin{equation*}
    (d\eta)(Z_1,\ldots, Z_{p+1}) = \sum_{i=1}^{p+1} (-1)^{i+1}
    (\nabla_{Z_i} \eta ) (Z_1, \ldots, \hat{Z_i}, \ldots, Z_{p+1}).
\end{equation*}
So, in terms of the family of functions we have defined above, using
Lemma \ref{nabla}, we see that
\begin{equation*}
    (d\eta)^{\alpha}_{i_1, \ldots, i_{p+1}} =\sum_{k=1}^{p+1}
    (-1)^{k+1} (X_{i_k} + \pi_1(X_{i_k}))
    \eta^{\alpha}_{i_1,\ldots, \hat{i_k}, \ldots, i_{p+1}}.
\end{equation*}

A simple argument (using Weil's integration by parts (equation
\ref{equation:weillifted}) and the fact that $\pi_1(X)$ is symmetric
with respect to the inner product on $F$ for $X\in \mathfrak{p}$
(Lemma \ref{F-admissible metric})) computes the adjoint of $d$ to be
as follows.

\begin{lem}
  $(\delta\eta)^{\alpha}_{i_1, \ldots, i_{p-1}} = -\sum_{k=1}^N
  (X_k - \pi_1(X_k))\eta^{\alpha}_{k,i_1,\ldots, i_{p-1}}.$
\end{lem}

Then we compute directly.
\begin{equation*}
\begin{split}
    (\delta d \eta)^{\alpha}_{i_1,\ldots,i_p} & =
    -\sum_{k=1}^N (X_k - \pi_1(X_k) ) (d\eta)^{\alpha}_{k, i_1,\ldots, i_p} \\
    & = -\sum_{k=1}^N (X_k -\pi_1(X_k))(X_k + \pi_1(X_k))\eta^{\alpha}_{i_1, \ldots,
    i_p}\\
    & \quad -\sum_{k=1}^N \sum_{j=1}^p (-1)^{j}(X_k -\pi_1(X_k))  (X_{i_j} + \pi_1(X_{i_j}))
    \eta^{\alpha}_{k, i_1, \ldots, \hat{i_j}, \ldots, i_p}.
\end{split}
\end{equation*}
\begin{equation*}
\begin{split}
    (d \delta \eta)^{\alpha}_{i_1, \ldots, i_p}  & = \sum_{j=1}^{p}
    (-1)^{j+1}
    (X_{i_j} + \pi_1(X_{i_j})) (d^*\eta)^{\alpha}_{i_1, \ldots, \hat{i_j}, \ldots,
    i_p} \\
    & =\sum_{j=1}^{p} (-1)^{j+1} (X_{i_j} + \pi_1(X_{i_j}) \left\{
    -\sum_{k=1}^N (X_k - \pi_1(X_k)) \eta^{\alpha}_{k,i_1, \ldots,
    \hat{i_j},\ldots, i_p} \right\} \\
    &= \sum_{k=1}^N \sum_{j=1}^p (-1)^{j} (X_{i_j}+ \pi_1(X_{i_j}))(X_k -
    \pi_1(X_k)) \eta^{\alpha}_{k, i_1 \ldots, \hat{i_j}, \ldots,
    i_p}.
\end{split}
\end{equation*}

We add these equations to compute the Laplacian.
\begin{equation}
\begin{split}
   (\triangle \eta)^{\alpha}_{i_1,\ldots, i_p} =& -\sum_{k=1}^N (X_k^2 -\pi_1(X_k)^2)
   \eta^{\alpha}_{i_1, \ldots, i_p} \\
     & + \sum_{k=1}^N \sum_{j=1}^p (-1)^j ([X_{i_j}, X_k] - \pi_1([X_{i_j},
     X_k])) \eta^{\alpha}_{k, i_1, \ldots, \hat{i_j}, \ldots, i_p},
\end{split}
\end{equation}
where we have used the equality $$(X_i + \pi_1(X_i))(X_k -
\pi_1(X_k)) -(X_k - \pi_1(X_k))(X_i + \pi_1(X_i)) = [X_i, X_k] -
\pi_1([X_i, X_k])$$ $p$ times to simplify the expression.

We may then proceed in the same way as Matsushima and Murakami do
\cite[page 385]{Matsushima-Murakami} to write the Laplacian in a
more algebraic form. Where Matsushima and Murakami use their
equation (4.10), we use the equivalent statement which follows
from propety P-IV of our lifted form $\eta^*$. This yields:
$$(\triangle \eta)^{\alpha}_{i_1,\ldots, i_p} = -C
\eta^{\alpha}_{i_1,\ldots,i_p} +
\pi_1(C)\eta^{\alpha}_{i_1,\ldots,i_p},$$ where $C = \sum_{k=1}^{N}
X_k^2 - \sum_{a=N+1}^n X_a^2$ is the Casimir operator of
$\mathfrak{g}$ and $\pi_1(C) =\sum_{k=1}^{N} \pi_1(X_k)^2 -
\sum_{a=N+1}^n \pi_1(X_a)^2$ is the Casimir operator of the
representation $\pi_1$. In terms of the $E=F\otimes\fh$-valued
functions $\eta_{i_1,\ldots,i_p}$ on $G$, we conclude that
\begin{equation}\label{equation:flatlaplacian}
    (\triangle \eta)_{i_1,\ldots,i_p} = -C \eta_{i_1,\ldots,i_p} + (\pi_1(C)\otimes
    \mathrm{Id}_{\mathcal{H}}) \eta_{i_1,\ldots,i_p}.
\end{equation}

\section{Obstructions to the approach via relative Lie algebra cohomology} \label{section:gkcohomology}

In this section, we outline briefly what would be required to
pursue an approach to our results via relative Lie algebra
cohomology (also called $(\mathfrak g, K)$-cohomology). We first
discuss direct approaches for connected groups, and then turn to
using Matsushima's formula for cocompact lattices.

\subsection{Connected Groups}
\label{subsection:gkcohomologydirect}  Given the fact, mentioned
following Theorem \ref{theorem:ft}, that vanishing of
$H^1(G,\pi{\otimes}\pi_2, F \otimes \fh)$ is known for $\pi_2$
irreducible or admissible, it is perhaps surprising that one
cannot deduce Theorem \ref{theorem:ft} from this. We describe the
scheme when $\pi_2$ is an infinite direct sum, the case of direct
integrals is essentially the same. If we have $(\pi_2,\fh) =
(\oplus_i \pi_2^i, \oplus \fh^i)$ then:
\begin{equation*}H^1(G, {\pi_1} {\otimes}
\pi_2, F{\otimes}\fh)= \frac{\oplus_i Z_1(G, \pi_1 \otimes \pi^i_2,
F{\otimes}\fh^i)}{\oplus_i B_1(G, \pi_1 \otimes \pi^i_2,
F{\otimes}\fh^i)},
\end{equation*}
where we know that the natural map $B_1(G, \pi_1 \otimes \pi^i_2,
F{\otimes}\fh^i) \rightarrow Z_1(G, \pi_1 \otimes \pi^i_2,
F{\otimes}\fh^i)$ has an inverse $p_i$ for every $i$. It is
straightforward to check that vanishing of the group
$H^1(D,{\pi_1}{\otimes}\pi_2, F{\otimes}\fh)$ is equivalent to a
uniform bound on the norms of all $p_i$.

The existing proofs that cohomology vanishes for $\pi_2$
admissible use quite involved homological arguments first to
translate the question to one of relative Lie algebra cohomology
and then to prove vanishing of relative Lie algebra cohomology. It
is not clear that the required bound can be produced by this
method.

One can try to bound $p_i$ a posteriori by direct computation,
using that $G$ has property $(T)$ and that $\pi_1$ is fixed. This
direction does not seem fruitful.

\subsection{Lattices}

We now discuss the possibility of a proof for cocompact lattices
using (an extension of) Matsushima's formula, a topic on which we
are more optimistic. Though it has a different guise here than
above, the problem is again one of uniformity of estimates. For this
section let $G$ be any semisimple Lie group with finite center. For
$\pi:\G{\rightarrow}GL(F)$ any finite dimensional representation
which almost extends to $G$, Matsushima's formula \cite{Matsushima}
says:
\begin{equation}
\label{equation:matsushima} H^k(\G,F)=\bigoplus n_{\pi}
H^k(\mathfrak g, \mathfrak k, V_{\pi,0}{\otimes}F),
\end{equation}
where the sum is over the irreducible $G$ representations $\pi$
appearing in $L^2(G/{\Gamma})$ and $V_{\pi,0}$ is the space of
$K$-finite vectors in the representation space $V_{\pi}$ for $\pi$.
The Hodge theorem and/or more formal variants are used in the course
of the proof.  In particular, the isomorphism in equation
\ref{equation:matsushima} above is deduced by computing cohomology
by computing harmonic forms.


In proofs of Matsushima's formula, there is a step at which there is
no isomorphism at the level of chain complexes, but only an
isomorphism between spaces of the harmonic forms in the two
complexes. In geometric language, the transition is between
computing cohomology using the complex of smooth forms and computing
cohomology using the complex of $L^2$ forms.  Since $L^2$ harmonic
forms are smooth, for finite dimensional vector bundles, the Hodge
theorem implies that one can compute cohomology by computing the
space of $L^2$ harmonic forms. The passage to $L^2$ forms is
important since it allows one to bring to bear various results from
representation theory in the computations. For a clear general
exposition of Matsushima's formula, see section $2$ of Schmid's
article \cite{Schmid-Advances}.

Since the Hodge theorem fails in general for infinite dimensional
vector bundles, it is no longer possible in our context to compute
cohomology using harmonic forms. However, we can produce a formula
along the lines of Matsushima's result which controls the allowable
harmonic forms in our more general bundle.

Let $G$ be a semisimple Lie group, $\G<G$ a cocompact lattice and
$\pi$ a finite dimensional by unitary representation of $\G$. To
avoid discussion of direct integrals, we assume that the unitary
representation $\pi_2$ is irreducible. Let ${\mathcal H}^k(\G,\pi,
F{\otimes}\fh)$ be the space of smooth harmonic $k$ forms in the
bundle $\G{\backslash}(G/K{\times}_{\pi}F{\otimes}\fh)$. Following
the proof of Matsushima's formula (with equation
\ref{equation:flatlaplacian} in place of Matsushima's computation of
the similar object), we obtain an isomorphism:
\begin{equation}
\label{equation:matsushimsaprime} {\mathcal H}^k(\G,\pi,
F{\otimes}\fh) = \bigoplus n_{\eta} H^k(\mathfrak g, \mathfrak k,
V_{\eta,0}{\otimes}F),
\end{equation}
where the sum is over the irreducible $G$ representations $\eta$
appearing in $I^G_{\G}(\pi_2)$ and $V_{\eta,0}$ is the space of
$K$-finite vectors in the representation space $V_{\eta}$ for
$\eta$. The expression above is correct because $H^k(\mathfrak g,
\mathfrak k, V_{\pi,0}{\otimes}F)$ can be computed as the cohomology
of a finite dimensional complex where a Hodge theorem is available
for purely formal reasons.

Combined with the theorem of Borel-Wallach, Zuckerman and Schmid
(precisely, the statement of \cite[Theorem 2.13]{Schmid-Advances}),
equation \ref{equation:matsushimsaprime} implies that ${\mathcal
H}^1(\G,\pi, F{\otimes}\fh)=0$ in the situation of Theorem
\ref{thm:1-vanishing}, so it is possible to control the part of the
cohomology which consists of harmonic forms.

It should be possible to give an alternate proof of Theorem
\ref{thm:1-vanishing} using an argument close to the one given by
Mok in \cite{Mok}.  We briefly sketch the idea to point out the
difficulties.  Given a bundle where an estimate of the type proven
in Theorem \ref{sufficient condition for vanishing} fails, Mok's
idea is to take a sequence of forms with $< \triangle f_i, f_i>
\rightarrow 0$ and $\|f_i\|_2=1$ and use a renormalizing and
limiting argument to construct another bundle over
$\G{\backslash}G/K$ where there are non-vanishing harmonic forms.
Mok does this for flat bundles associated to unitary
representations, and shows that the bundle constructed in the limit
has the same form.  In order to have a result applicable in our
setting, we would need to repeat Mok's argument with a flat bundles
associated to a finite dimensional by unitary representation and
prove that the limiting representation is finite dimensional by
unitary.  We believe this is possible and hope to pursue it in a
later paper, where we would also give a detailed account of equation
\ref{equation:matsushimsaprime} and it's generalizations. These two
results would imply better vanishing results for higher degree
cohomology than those in \S \ref{subsection:higherdegrees}. In fact,
in combination with results of Vogan-Zuckerman, we expect this will
yield optimal vanishing theorems for higher degree cohomology of
cocompact lattices \cite{Vogan-Zuckerman}.  However, due to the
failure of Lemma \ref{lemma:restriction} for $H^k$ when $k>1$, this
will not yield vanishing theorems for the ambient group.

\bibliographystyle{amsplain}
\bibliography{CSRproject}
\end{document}